\documentclass[12pt]{amsart} 
\usepackage{amsthm}  
\usepackage{epsf}  
\usepackage{amssymb}  
\usepackage{latexsym}  
\usepackage{amsmath}  
\newtheorem{thm}{Theorem}[section] 
\newtheorem{defn}[thm]{Definition}  
  
\newtheorem{cor}[thm]{Corollary}  
\newtheorem{lemma}[thm]{Lemma}  
  
\newtheorem{prop}[thm]{Proposition}  
\newtheorem{rem}[thm]{Remark}  

\newcommand{\bde}{\begin{defn}}
\newcommand{\ede}{\end{defn}}
\newcommand{\bpr}{\begin{prop}}
\newcommand{\epr}{\end{prop}}
\newcommand{\bthm}{\begin{thm}}
\newcommand{\ethm}{\end{thm}}
\newcommand{\ben}{\begin{enumerate}}
\newcommand{\een}{\end{enumerate}}
\newcommand{\ble}{\begin{lemma}}
\newcommand{\ele}{\end{lemma}}

\def\R{{\Bbb R}}  
  
\def\C{{\Bbb C}}

\def\PP{{\Bbb P}}

\def\cpt{\Bbb C \Bbb P^2}  
\def\cpn{\Bbb C \Bbb P^N}

\def \ta{\tau}  
\def \ta1{\tau_1}

\def \be{\beta}  
  
\def \de{\delta}  

\def \ga{\gamma}  
  
\def \Ga{\Gamma}

\def \sbs{\subset}

\def \cB{{\mathcal B}}  
  
\def \cR{{\mathcal R}}

\def \zovera {  
    \mathop{\lower 10pt \hbox{${\buildrel{\displaystyle\bar{z}} \over {\scriptstyle{(a)}}} $}}  
    {\lower 4pt \hbox{${\scriptstyle{ij}}$}} 
} 

\newcommand\figs[1]{#1}

\newcommand\FIGUREy[4][]{{\begin{figure}[!h]\epsfysize=#3 \centerline{\epsfbox{\figs{#2}}}
\caption{#1}\label{#4}\end{figure}}}

\newcommand\begintable[1][] {{}}

\newif\ifXY 
\XYtrue     
\newif\ifbigmatrices  
\bigmatricestrue 
  
 \ifXY  
 \usepackage{xy}  
 \fi  
 \ifXY  
 \xyoption{all}  
 \fi  
  
\begin{document}  
  
\title[Local braid monodromies for conic-line arrangements] {Local braid monodromies 
and local fundamental groups of tangented conic-line arrangements}  
  
\author[Meirav Amram, David Garber and Mina Teicher]{Meirav Amram$^1$, David Garber$^2$ and Mina Teicher}  
  
\stepcounter{footnote} 
\footnotetext{Partially supported by the Minerva Foundation of Germany, the DAAD 
fellowship (Germany), EAGER (EU Network, HPRN-CT-2009-00099) and 
the Golda Meir postdoctoral fellowship.}
\stepcounter{footnote} 
\footnotetext{Partially supported by the Lady Davis postdoctoral fellowship.} 

\address{Meirav Amram, Einstein Institute of Mathematics, 
The Hebrew University, Givat Ram, 91904 Jerusalem, Israel}
\email{ameirav@math.huji.ac.il}

\address{David Garber, Einstein Institute of Mathematics, 
The Hebrew University, Givat Ram, 91904 Jerusalem, Israel}

\email{garber@math.huji.ac.il}

\address{Mina Teicher, Department of Mathematics, Bar-Ilan University, 52900
Ramat Gan, Israel}  
  
\email{teicher@macs.biu.ac.il}

  
 \renewcommand{\subjclassname}{%
       \textup{2000} Mathematics Subject Classification}  
  
  
 \date{\today}

\begin{abstract}  
This paper is the first part of a series of three papers about the fundamental
groups of conic-line arrangements consist of two tangented conics 
and up to two additional tangented lines. 

In this part, we compute the local braid monodromies and the local fundamental 
groups of the singularities which appear in such arrangements.
\end{abstract}  
  
\maketitle  
  
\section{Introduction}\label{intro}  
  
There is an extensive research concerning the   
topology of line arrangements and curves, and the fundamental groups
of their complements, see \cite{AmTe}, \cite{AmTeUl}, \cite{Ar},
\cite{De}, \cite{Fu}, \cite{GaTe}, \cite{Li1}, \cite{Li2},
\cite{MoTe}, \cite{MoTe5}-\cite{Sa}, and others.   

The study of fundamental groups of complements of algebraic curves   
contributes to other areas of algebraic geometry.   

In \cite{Ku},\cite{KuTe} and \cite{Te2}, we presented invariants of branch curves 
arising from generic projections of surfaces in $\cpn$ to $\cpt$ as 
distinguishing invariants between different connected components of 
the moduli space of surfaces of general type. 
Moreover, in \cite{Te}, we presented fundamental groups of complements of branch 
curves as new invariants for surfaces.

This paper is the first part of a series of three papers which 
study the fundamental group of some types of conic-line arrangements, which may appear 
as a branch curve of a generic projection of a surface 
of general type in $\cpn$ to $\cpt$ (see for example 
\cite{Hi},\cite{MRT},\cite{Te},\cite{Te2} and other works in progress 
\cite{ACMT},\cite{CMT}). The other two parts are \cite{AmGaTeII}
and \cite{AmGaTeIII}.
  
In particular, we study singularities which appear in  conic-line
arrangements with two tangented conics and up to two additional 
tangented lines.

\medskip
   
The paper is organized as follows. Section \ref{bm} surveys briefly the notion
of braid monodromy and the way to compute it, 
and it includes the presentation of the van Kampen Theorem which computes
the induced relations from the braid monodromy of the singular point.   
In Section \ref{2_comps} we compute the local 
braid monodromy and the relations induced by the simple tangency point, 
created by two tangented conics. Section  \ref{3_comps} deals with the local braid 
monodromy and the relations induced by singular points with three components. 
Section \ref{4_comps} covers some singular points with four components.

\section{The computation of local braid monodromy and its induced relations}\label{bm}

In this section, we briefly survey the notion of braid monodromy of a 
singular point of a curve, and the way how we deduce its 
corresponding relations in the fundamental group.  

We start with the definition of $l( \ga )$ and a {\it g-base}, 
which will be needed later. Let $K$ be a finite set in a disk $D \sbs \R^2$.

\bde \label{l_si}
Let $D$ be a disk. Let $w_i, \ i=1, \cdots , n$, be  disjoint disks in ${\rm Int}(D)$.
Let $u \in \partial D$.
Let $\ga$ be a simple path connecting $u$ with one of the $w_i$s, 
say $w_{i_0}$, and does not meet any other $w_j$. 
We assign to $\ga$ a loop $l( \ga )$ (actually an element of $\pi _1(D-K,u)$)
as follows: let $c$ be a simple loop equal to the (oriented) boundary of a 
small neighbourhood $V$ of $w_{i_0}$ chosen such that $\ga ' = \ga - V \cap \ga$ is
a simple path (see Figure \ref{gbase-fig}). \\
 Then: $l( \ga ) = \ga ' \cup c \cup (\ga ')^{-1}$ (we will not 
distinguish between $l( \ga )$ and its representative in $\pi _1(D-K,u)$).
\ede

\FIGUREy[$l(\ga)$]{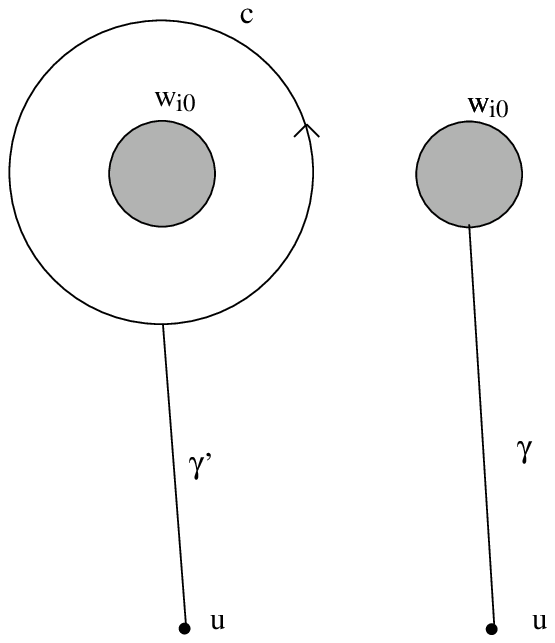}{4cm}{gbase-fig}

\bde \label{g_base} 
Let $D$ be a disk, $K \sbs D$, $| K| < \infty$. Let $u \in D-K$.  
A set of simple paths $\{ \ga _i \}$ is a {\bf bush} in $(D,K,u)$, if 
$\forall i,j,\ \ga _i \cap \ga _j = u; \ \forall i,\ \ga _i \cap K$ 
contains only one point, and $\ga _i$ are ordered counterclockwise around u. 
Let $\Ga _i = l(\ga _i) \in \pi _1(D-K,u)$ be a loop around $K \cap \ga _i$ 
determined by $\ga _i$. $\{ \Ga _i \}$ is called a {\bf g-base} of $\pi _1(D-K,u)$. 
\ede

An example of a g-base is given in Figure \ref{page14}. 

\FIGUREy[A g-base]{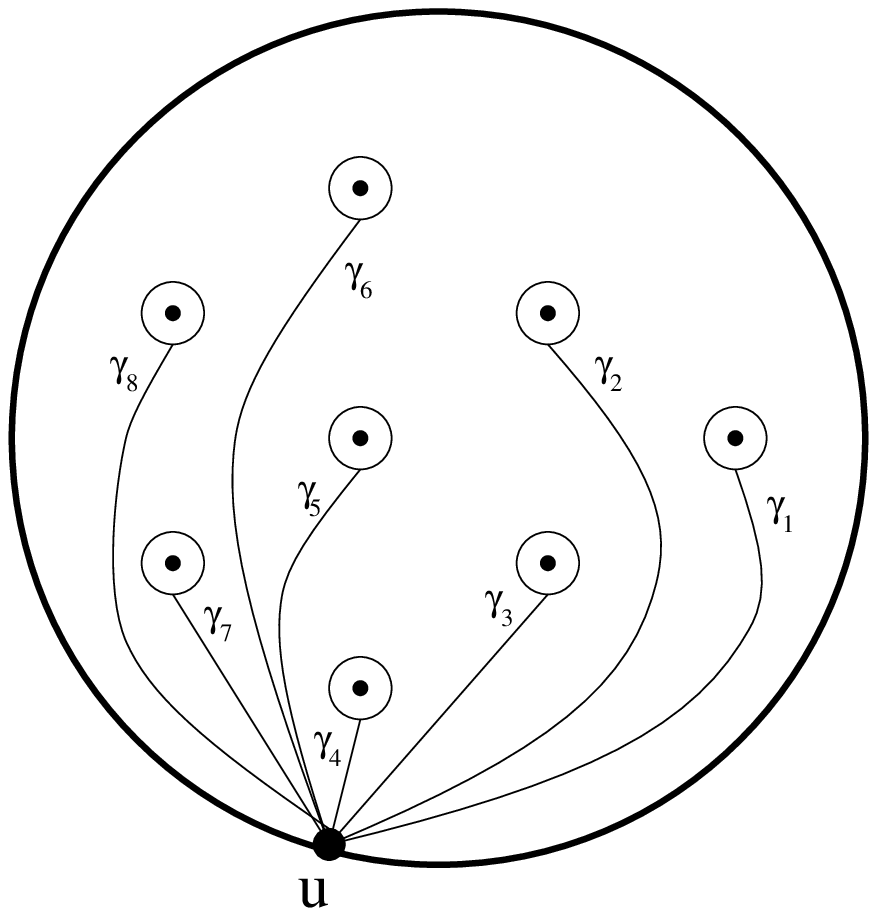}{5cm}{page14}

Here we define the {\it standard g-base} which will be used throughout the paper:

\bde \label{std_g_base} 
Let $D$ be the unit disk, and let $K\sbs D$ be a finite set of points on 
the real axis. Let $u = (-1,0)$.  
Let $\{ \ga _i \}$ be the straight segment which connects $u$ to the $i$th
point. Let $\Ga _i = l(\ga _i) \in \pi _1(D-K,u)$. $\{ \Ga _i \}$ 
is called the {\bf standard g-base} of $\pi _1(D-K,u)$ (see Figure \ref{std_gbase}).
 
\FIGUREy[A standard g-base]{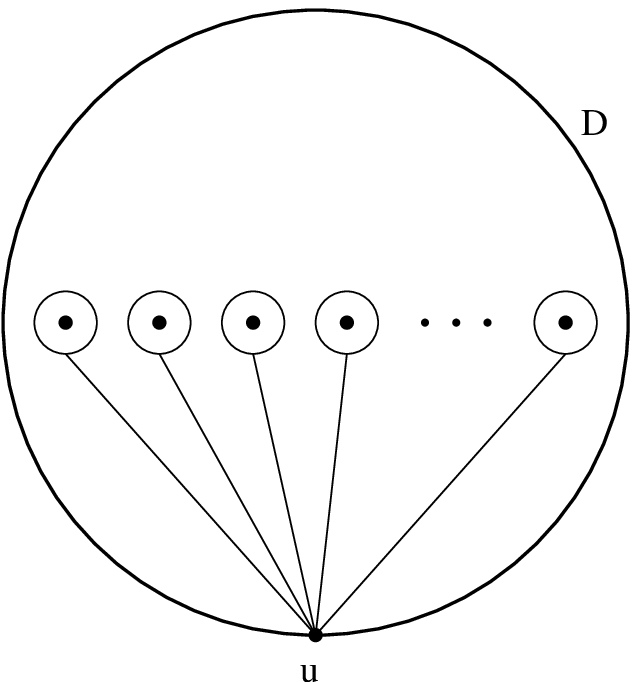}{4cm}{std_gbase}
\ede

\subsection{The braid monodromy and the way to compute it}
Let $D$ be again a closed disk in $\R^2$, $K \sbs D$ a finite set, $u \in\partial D$.
In such a case, we can define the {\it braid group} $B_n[D,K]\ (n= |K|)$:

\bde \label{BG-topology} 
Let $\cB$ be the group of all diffeomorphisms $\be$ of $D$ such that $\be(K)=K$,
$\be|_{\partial D} = {\rm Id}|_{\partial D}$. Such a diffeomorphism acts naturally on
$\pi_1(D-K,u)$ (since such a diffeomorphism acts on paths in general, and on loops in
particular). We say that two such diffeomorphisms
are equivalent if they define the same automorphism on $\pi_1(D-K,u)$.
The quotient of $\cB$ by this equivalence relation is called the 
{\bf braid group $B_n[D,K]$}. 
\ede

An alternative equivalent definition for the braid group was given by 
Artin \cite{Artin}, by means of generators and relations:

\bde[Artin] \label{BG-artin} 
The {\bf braid group $B_n$} is the group generated by $\{\sigma _1,...,\sigma _{n-1}\}$
subject to the relations: 
\ben
\item \ $\sigma _i\sigma _j=\sigma _j\sigma _i$ where $|i-j| \geq 2$. 
\item \ $\sigma _i\sigma _{i+1}\sigma _i=\sigma _{i+1}\sigma _i\sigma _{i+1}$  
for all $i=1,\dots,n-2$.
\een
\ede

Another way to look at braids is via {\it motions} of $K$.
 
\bde 
Let $K'= \{ a'_1, a'_2, \cdots ,a'_n \}, \ K= \{ a_1,a_2, \cdots ,a_n \}$. 
{\bf A motion of K' to K in D} is a set of $n$ continuous functions 
$m_i:[0,1] \to D, \ i=1, \cdots ,n$, such that:  
\begin{enumerate}
\item $\forall i,\  m_i(0) = a'_i, \ m_i(1) = a_i$. 
\item $\forall i \not = j, \ m_i(t) \not = m_j(t) \ \forall t \in [0,1]$. 
\end{enumerate}
\ede
 
According to the following proposition, we can define a family of
diffeomorphisms induced from the motion (with the assumption that $K=K'$).
  
\bpr
Given a motion $\cR$, there exists a continuous family of diffeomorphisms 
$D_{\cR ,t} : D \to D, \ t \in [0,1]$, such that:
\begin{enumerate}
\item $D_{\cR ,t}|_{\partial D} = {\rm Id}|_{\partial D}$. 
\item $\forall t,i,\  D_{\cR ,t}(a'_i)=m_i(t)$.
\end{enumerate}
\epr
 
\bde \label{b_m} 
When $K=K'$, the braid induced from a motion $\cR$, denoted by $b_\cR$, 
is the braid defined by the diffeomorphism $D_{\cR ,1}$.
\ede

Before introducing the definition of {\it braid monodromy}, we have to make
some more constructions. From now on, we will work in $\C^2$. Let $E$ (resp. $D$) 
be a closed disk on the $x$-axis (resp. $y$-axis), 
and let $S$ be a part of an algebraic
curve in $\C^2$ located in $E \times D$. Let $\pi _1 : E \times D \to E$ and
$\pi _2 : E \times D \to D$ be the canonical projections, and let 
$\pi = \pi _1 |_S: S \to E$. Assume $\pi$ is a proper map, and $\deg \pi = n$. 
Let $ N = \{ x \in E \ |\ | \pi ^ {-1}(x)| < n \} $, and assume 
$ N \cap \partial E= \emptyset $. Now choose $ M \in \partial E$ 
and let $K = K(M) = \pi ^ {-1} (M)$. By the assumption that 
$\deg \pi = n \ \ (\Rightarrow |K|=n)$,  we can write: 
$K = \{ a_1,a_2, \cdots ,a_n \} $.
Under these constructions, from each loop in $E-N$, we can define a braid in
$B_n [M \times D, K]$ in the following way: 
\begin{itemize}
\item[(1)] Because $ \deg \pi = n$, we can lift any loop in $E-N$ with 
a base point $M$ to a system of $n$ paths in $ (E-N) \times D $ which 
start and end at $ \{ a_1,a_2, \cdots ,a_n \} $.
\item[(2)] Project this system into $D$ (by $\pi _2$), to get $n$ paths in $D$ 
which start and end at the image of $K$ in $D$ (under $\pi_2$). These paths 
actually form a motion.
\item[(3)] Induce a braid from this motion, as defined in Definition \ref{b_m}.
\end{itemize}

In this way, we can match a braid to each loop. 
Therefore, we get a map $ \varphi : \pi _1(E-N,M) \to B_n [M \times D, K]$,
which is also a group homomorphism, called the
{\bf braid monodromy of $S$ with respect to $E \times D,\pi _1,M$}.      

Hence, for computing the braid monodromy, one should perform the following steps: 
\begin{enumerate}
\item Define a loop in the $x$-axis which circumscribes the $x$-projection of 
the singular point. 
\item Compute the points of the curve $S$ in each fiber over the loop.
\item Determine the behavior of the points in the fibers which correspond 
to the loop. This gives the braid monodromy.   
\end{enumerate}  

Some examples of computing braid monodromies of singular points 
can be found in \cite{MoTe1}.

\subsection{The Lefschetz diffeomorphism}

For computing the braid monodromy of a plane curve, we use the 
Moishezon-Teicher algorithm (see \cite{MoTe1}). 
The idea of this algorithm is to start with some braid which 
is associated to the local behavior near the  
singular point, and then to apply on it some actions induced by 
the singularities to the right side of that singular point. 
This actions are called {\it Lefschetz diffeomorphisms} which are induced 
by the singular points.  

\begin{figure}[h]
\epsfysize=2.5cm  
\centerline{\epsfbox{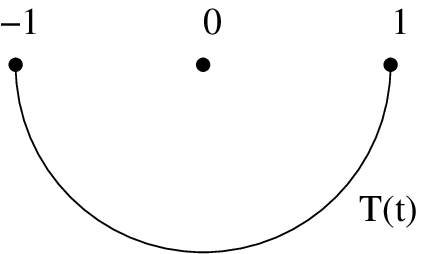}}  
\caption{The path $T(t)$}\label{path}  
\end{figure}

The Lefschetz diffeomorphism can be defined in general (see \cite{MoTe1}), 
but here we define it according to our settings. Let $P$ be a singular point, 
and assumed that $P$ is located at the origin. Let $T(t), 0 \leq t \leq 1$, 
be the path  which connects the point $T(0)=-1$ to the point $T(1)=1$ 
along the negative part of the unit circle (see Figure \ref{path}).   

There exists a continuous family of diffeomorphisms 
$\psi _{(t)}: \{ -1 \} \times D \to T(t) \times D, \ \forall t \in [0,1]$, 
such that: 
\begin{enumerate}
\item $\psi _{(0)} = {\rm Id}|_{ \{ -1 \} \times D}$. 
\item $\forall t \in [0,1], \ \psi _{(t)}(\pi^{-1} _1( \{ -1 \} ) \cap S) =
                            \pi^{-1} _1(T(t)) \cap S$.
\item $\forall y \in \partial D, \ \psi _{(t)}(\{ -1 \},y) = (T(t),y)$.
\end{enumerate}

\bigskip

In this situation, we can define the {\it Lefschetz diffeomorphism induced by the singular point $P$}:

\bde 
The {\bf Lefschetz diffeomorphism induced by the singular point $P$} is the diffeomorphism
$$\psi _T = \psi _{(1)}: \{ -1 \} \times D\  \tilde\to \ \{ 1 \} \times D $$
\ede

One can see that the computation of the Lefschetz diffeomorphism is similar 
to the computation of the braid monodromy, in the following sense: in order to compute 
the braid monodromy, we have to follow the movement of the points of the curve 
in the fibers over the whole loop, whence for computing the Lefschetz diffeomorphism, 
this following should be done only for a half of the loop. 

\subsection{The classical van Kampen Theorem}\label{vK-thm} 

The van Kampen Theorem gives a finite presentation (by generators and relations) of the 
fundamental group of complements of curves from the braid monodromy. 

\medskip

Let $S$ be an algebraic curve in $\C^2$ ($\deg S =n$). 
Let $\pi = \pi_1 : \C ^2 \to \C$
be the canonical projection on the first coordinate.   
Let $\C _x = \pi ^{-1}(x)$, and now define: $K_x=\C _x \cap S$ (By assumption
${\rm deg}S=n$, we have that $| K_x | \le n$). 

Let $N=\{x\ |\ | K_x | < n \} $.
Choose now $u \in \C$, $u$ real, such that $x \ll u, \ \forall x \in N$, and define:
$B_n = B_n [ \C _u, \C _u \cap S ]$. Let $\varphi _u : \pi _1(\C -N, u) \to B_n$
be the braid monodromy of S with respect to $\pi, u$.
Also choose  $u_0 \in \C _u, \ u_0 \not\in S$, $u_0$ below the real line far
enough such that $B_n$ does not move $u_0$.
It is known that the group $\pi _1(\C _u - S, u_0)$ is free.
There exists an epimorphism
$\pi _1(\C _u - S, u_0) \to \pi _1(\C ^2 - S, u_0)$, so a set of generators 
for $\pi _1(\C _u - S, u_0)$ determines a set of generators for      
$\pi _1(\C ^2 - S, u_0)$. 

\medskip

In this situation, the classical van Kampen Theorem says:

\bthm[Classical van Kampen's Theorem] \label{VK1} 
Let $S$ be an algebraic curve, $u,u_0,\varphi _u$ 
defined as above. Let $\{ \de _i \}$ be a g-base of $\pi _1(\C -N,u)$. 
Let $\{ \Ga _j \ | \ 1 \le j \le n \}\ (n=\deg S)$ be a g-base for
$\pi _1 (\C _u -S,u_0)$. \\
Then, $ \pi _1(\C ^2 -S, u_0)$ is generated
by the images of $\Ga _j$ in $\pi _1(\C ^2 -S, u_0)$ and we get a 
complete set of relations
from those induced from
$$(\varphi _u(\de _i))(\Ga _j) = \Ga _j, \qquad \forall i \forall j$$
\ethm

Hence, for computing the relations induced from a singular point, one should 
follow these steps (see also \cite{MoTe1}-\cite{MoTe2}):
\begin{enumerate}
\item Start from a standard g-base.
\item Apply on the g-base the corresponding braid (given by the braid monodromy). 
\item Present any loop as a composition of the generators of the standard g-base.
\item The induced relations are the $n$ equations: for all $j$, $\Ga_j$ is equal 
      to the presentation of the $j$th loop.
\end{enumerate}

One should notice that all the computations here are {\it local}, and when we reach
the global computations, we need to replace the generators appear in these relations
by a word (which is a conjugation of a standard generator by some 
standard generators).

\section{A singular point with two components}\label{2_comps}

In this section, we compute the local braid monodromy of a tangency point of 
two conics.

\begin{figure}[h]
\epsfysize=4.5cm  
\centerline{\epsfbox{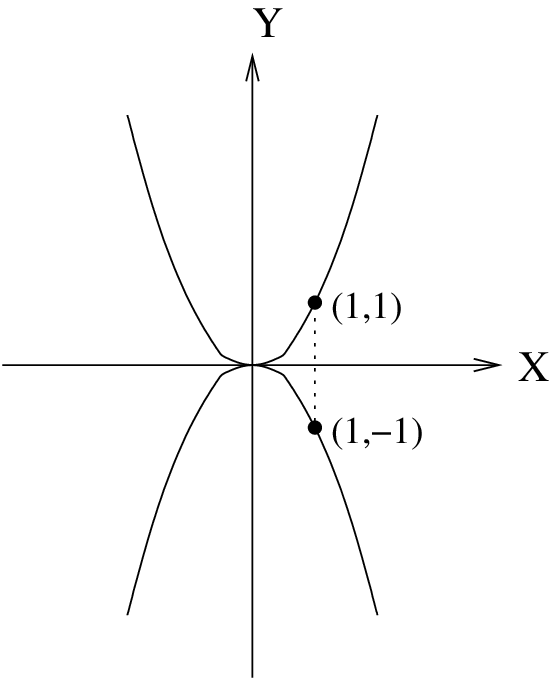}}  
\caption{The singularity of $(y+x^2)(y-x^2) = 0$ at $(0,0)$}\label{2TanCon}  
\end{figure}

The tangency point between two conics can be presented locally by the 
equation: $(y+x^2)(y-x^2) = 0$ (see Figure \ref{2TanCon}).  

The singular point is $(0,0)$, and the points of  
the curve in the fiber over $x=1$ are $y=1$ and $y=-1$. 
Now, let us take a loop in $x$-axis based at $x=1$ and circumscribing $0$:   
$\alpha (t) =e^{2 \pi it}$, where  $0 \leq t \leq 1$.  Lifting it to
the curve, we get two paths on the curve:  
$\alpha _1(t) = (e^{2 \pi it}, e^{4 \pi it})$ and   
$\alpha _2(t) = (e^{2 \pi it}, -e^{4 \pi it})$ where  $0 \leq t \leq 1$.  
Projecting the paths to the $y$-axis, we get $y_1(t) = e^{4 \pi it}$ and  
 $y_2(t)= -e^{4 \pi it}$ where $0 \leq t \leq 1$.  
By substituting $t=1$, the local braid monodromy of this point 
is two counterclockwise full-twists. Hence, we have:

\begin{thm}
The local braid monodromy of the tangency point of two conics is two counterclockwise
full-twists (similar to the tangency point of a conic and a line). Hence, we have 
the same induced relations: 
$$(x_1 x_2)^2 = (x_2 x_1)^2$$
where $\{ x_1,x_2 \}$ are the generators of the standard g-base.
\end{thm}

Now, we compute the corresponding Lefschetz diffeomorphism: 
\begin{cor}
The Lefschetz diffeomorphism of a tangency point between two tangented conics 
is a counterclockwise full-twist (similar to the usual tangency point).
\end{cor}

The reason for this corollary is since the Lefschetz diffeomorphism is obtained
by going on the loop mentioned above from $t=\frac{1}{2}$ 
to $t=1$. Along this interval, the two points in the fiber make one counterclockwise 
full-twist, and hence the Lefschetz diffeomorphism of this singular point
is a full-twist.

\section{Singular points with three components}\label{3_comps}

In this section, we will compute the local braid monodromy of four types 
of singular points consist of three components (where at most two of them 
are conics which are tangented). Then, we will compute 
the relations induced from these singular points by the classical van Kampen Theorem
(see \cite{vK} and Section \ref{vK-thm} here). 

\subsection{First type}\label{3comps_type1}
 
\begin{figure}[h]
\epsfysize=2.5cm  
\centerline{\epsfbox{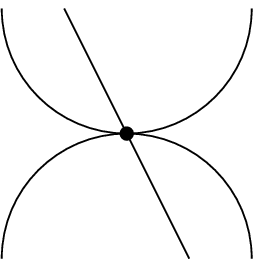}}  
\caption{The singularity of $(2x+y)(y+x^2)(y-x^2) = 0$ at $(0,0)$}\label{type_b}  
\end{figure}

The local equation of the singularity of the first type is 
$(2x+y)(y+x^2)(y-x^2) = 0$ (see Figure \ref{type_b}). 

For computing the braid monodromy, we take a loop around $x=0$ and look what happens
to the points of the curve $(2x+y)(y+x^2)(y-x^2) = 0$ in the fibers above this loop. 
 
Let $x=e^{2\pi i t}$ where $0 \leq t \leq 1$. For $t=0$ we have that $x=1$ and 
the points of the curve in the fiber over $x=1$ are $y=1,-1$ and $y=-2$. 
For $t=\frac{1}{2}$, we have $x=-1$, and the points of the curve 
in the fiber over $x=-1$ are $y=-1,1$ and $y=2$.  The point $y=1$
in the fiber $x=-1$ corresponds to the point $y=1$ in the fiber $x=1$, and similarly
the point $y=-1$ in the fiber $x=-1$ corresponds to the point $y=-1$ 
in the fiber $x=1$. The point $y=2$ in the fiber $x=-1$ corresponds to the point $y=-2$ 
in the fiber $x=1$. Hence, we get that from $t=0$ to $t=\frac{1}{2}$, the points $y=1$
and $y=-1$ do a counterclockwise full-twist and the point $y=-2$ does 
a counterclockwise half-twist around the points $y=1,-1$. 
If we continue to $t=1$, the points $y=1$
and $y=-1$ do two counterclockwise full-twists and the point $y=-2$ does 
a counterclockwise full-twist around the points $y=1,-1$ together. 

\begin{figure}[h]
\epsfysize=8cm  
\centerline{\epsfbox{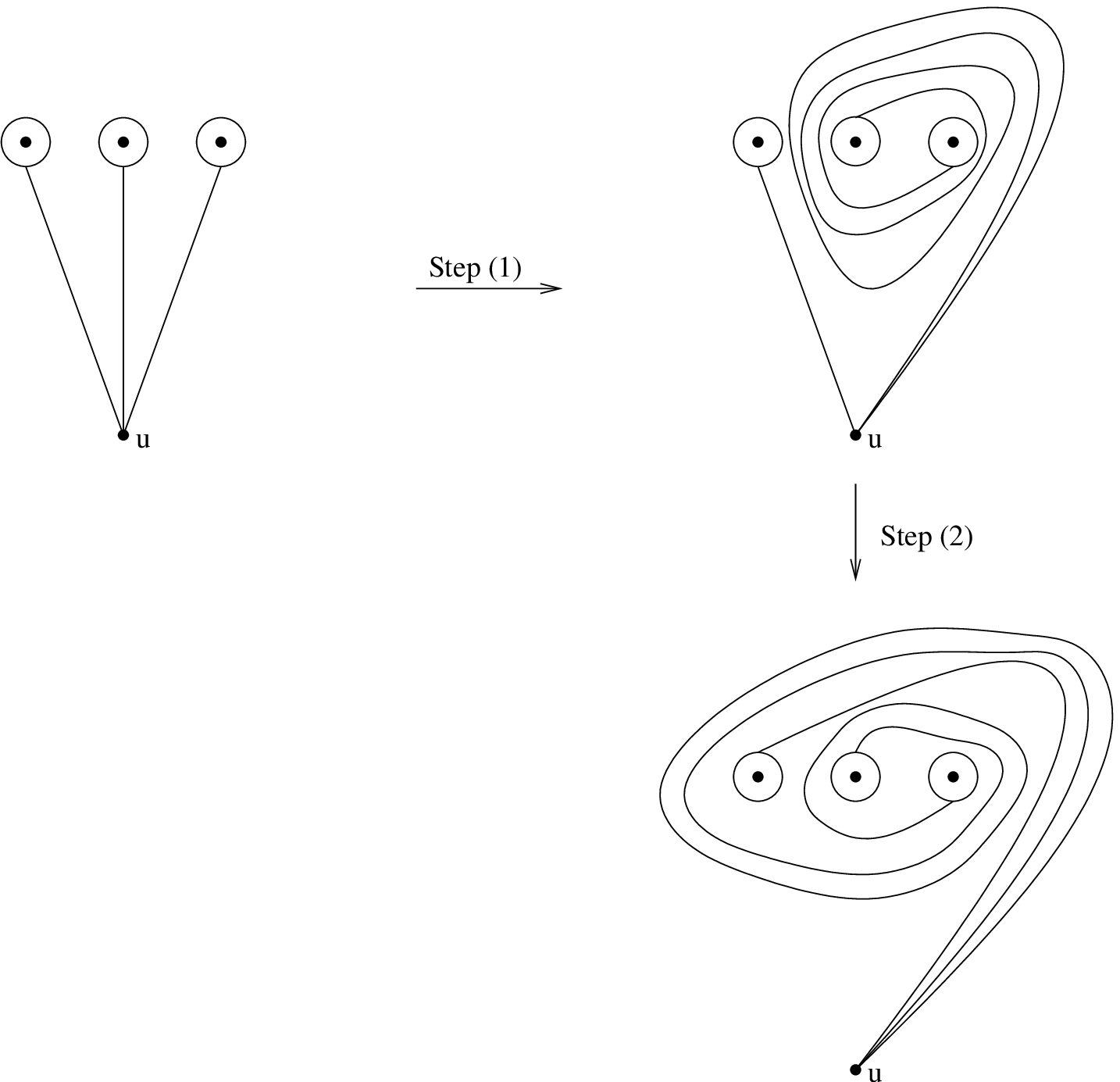}}  
\caption{The g-base obtained from the standard g-base by the action of the local braid monodromy}\label{gbase_type_b}  
\end{figure}

Now, we want to compute the induced relations from the braid monodromy 
of the singular point. According to van 
Kampen's Theorem, one should compute the g-base obtained by applying the action induced
by the local braid monodromy of the singular point on the standard g-base. 
This is shown in Figure \ref{gbase_type_b}: in the first step, we perform two 
counterclockwise full-twists of points $2$ and $3$ which represent 
$y=-1$ and $y=1$ respectively. In the second step, 
we perform a counterclockwise full-twist of point $1$ (which represents 
the point $y=-2$) around the points $2$ and $3$.

Now, by van Kampen's Theorem, we get the following induced relations from 
the new g-base (where $\{ x_1,x_2,x_3 \}$ are the generators of the standard g-base):
\begin{enumerate}
\item $x_1 = x_3 x_2 x_1 x_2^{-1} x_3^{-1}$  
\item $x_2 = x_3 x_2 x_1 x_3 x_2 x_3^{-1} x_1^{-1} x_2^{-1} x_3^{-1}$
\item $x_3 = x_3 x_2 x_1 x_3 x_2 x_3 x_2^{-1} x_3^{-1} x_1^{-1} x_2^{-1} x_3^{-1}$
\end{enumerate}

From Relation (1), we get that $x_1 x_3 x_2 = x_3 x_2 x_1$. Relation (2) becomes
$x_3 x_2 x_1 x_3 x_2 = x_2 x_3 x_2 x_1 x_3$. 

Now, we show that Relation (3) is not needed. Simplifying Relation (3) yields
$x_2 x_1 x_3 x_2 x_3 = x_3 x_2 x_1 x_3 x_2$. By Relation (1), 
this relation is equal to Relation (2), and hence Relation (3) is redundant.   

To summarize, we proved the following:

\begin{thm}
The local braid monodromy of the singularity presented locally by the equation 
$(2x+y)(y+x^2)(y-x^2)=0$ is: two points (correspond to $y=1$ and $y=-1$) do two 
counterclockwise full-twists and the third point (corresponds to $y=-2$) do 
a counterclockwise full-twist around them. 

The induced relations of this point are: 
$$x_1 x_3 x_2 = x_3 x_2 x_1 \qquad ; \qquad x_3 x_2 x_1 x_3 x_2 = x_2 x_3 x_2 x_1 x_3$$
where $\{ x_1,x_2,x_3 \}$ are the generators of the standard g-base.
\end{thm}

Since the Lefschetz diffeomorphism of the singular point is obtained by computing
the action only on half of the unit circle (from $t=\frac{1}{2}$ to $t=1$), 
we have the following corollary:

\begin{cor}
The Lefschetz diffeomorphism of the singular point presented locally by
$(2x+y)(y+x^2)(y-x^2) = 0$ is a counterclockwise full-twist of the points
correspond to $y=1$ and $y=-1$, and a half-twist of the point corresponds to $y=2$ 
(which becomes $y=-2$) around them.
\end{cor}

\begin{rem}
One can easily see that the singular point presented locally by the equation 
$y(2x+y)(y+x^2)=0$ (see Figure \ref{type_g2}) has the same braid monodromy 
and the same induced 
relations as the singular point which we have dealt with in this section, 
since these two singular points are locally the same (from the topological point 
of view).

\begin{figure}[h]
\epsfysize=2.5cm  
\centerline{\epsfbox{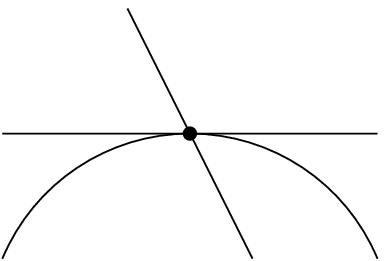}}  
\caption{The singularity of $y(2x+y)(y+x^2) = 0$ at $(0,0)$}\label{type_g2}  
\end{figure}  
  
\end{rem}

\subsection{Second type}\label{3comps_type2}

\begin{figure}[h]
\epsfysize=2.5cm  
\centerline{\epsfbox{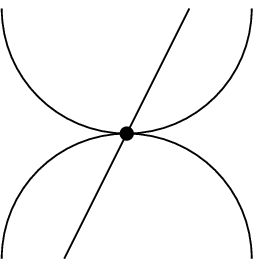}}  
\caption{The singularity of $(2x-y)(y+x^2)(y-x^2) = 0$ at $(0,0)$}\label{type_c}  
\end{figure}

The local equation of the singularity of the second type is 
$(2x-y)(y+x^2)(y-x^2) = 0$  (see Figure \ref{type_c}). 

Applying almost the same computations we have performed for the first type, 
we get the following result:

\begin{thm}
The local braid monodromy of the singularity presented locally by the equation 
$(2x-y)(y+x^2)(y-x^2)=0$ is: two points (correspond to $y=1$ and $y=-1$) do two 
counterclockwise full-twists and the third point (corresponds to $y=2$) do 
a counterclockwise full-twist around them. 

The corresponding induced relations are:    
$$x_3 x_2 x_1 = x_2 x_1 x_3 \qquad ; \qquad x_3 x_2 x_1 x_2 x_1 = x_1 x_3 x_2 x_1 x_2$$
where $\{ x_1,x_2,x_3 \}$ are the generators of the standard g-base.
\end{thm}

The Lefschetz diffeomorphism of the singular point 
is a counterclockwise full-twist of the points
correspond to $y=1$ and $y=-1$, and a half-twist of the point corresponds to $y=-2$ 
(which becomes $y=2$) around them.

\begin{rem}

As in the previous section, one can easily see that the singular point presented 
locally by the equation $y(2x-y)(y+x^2)=0$ (see Figure \ref{type_g}) 
has the same braid monodromy and the same induced 
relations as the singular point which we have dealt with in this section, 
since these two singular points are locally the same.

\begin{figure}[h]
\epsfysize=2.5cm  
\centerline{\epsfbox{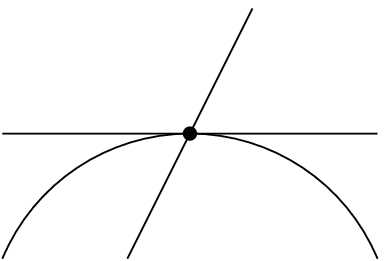}}  
\caption{The singularity of $y(2x-y)(y+x^2) = 0$ at $(0,0)$}\label{type_g}  
\end{figure}  
  
\end{rem}

\subsection{Third type}\label{3comps_type3}

\begin{figure}[h]
\epsfysize=2.5cm  
\centerline{\epsfbox{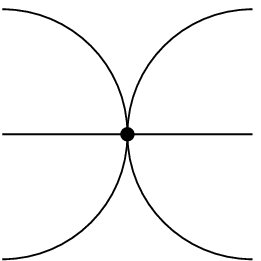}}  
\caption{The singularity of $y(y^2+x)(y^2-x) = 0$ at $(0,0)$}\label{type_h}  
\end{figure}

The local equation of the singularity of the third type is $y(y^2+x)(y^2-x) = 0$ 
(see Figure \ref{type_h}). One should notice that there is a major difference between
this type of singularity and the two previous types: In this singularity there are
two ``hidden'' branch points. That is, at any fiber one has three real points 
and two complex points (i.e. complex level $2$), and in each side of the 
singularity, the complex points belong to a different conic (by the 
singularity, two real points become complex and two complex points become real).   

For computing the braid monodromy, we take a loop around $x=0$ and look what happens
to the points of the curve $y(y^2+x)(y^2-x) = 0$ in the fibers over this loop. 
 
Let $x=e^{2\pi i t}$ where $0 \leq t \leq 1$. For $t=0$ we have that $x=1$ and 
the points of the curve in the fiber over $x=1$ are $y=-1,0,1,i$ and $y=-i$. 
For $t=\frac{1}{2}$, we have $x=-1$, and we get again that the points of the curve 
in the fiber over $x=-1$ are $y=-1,0,1,i$ and $y=-i$.  By a careful checking, 
one can see that the point $y=1$
in the fiber $x=-1$ corresponds to the point $-i$ in the fiber $x=1$, and similarly
all the points except for $y=0$ in the fiber $x=-1$ 
made a $90^{\circ}$ rotation counterclockwise around $y=0 \in \C$ 
from their corresponding points in the fiber $x=1$. The point $y=0$ remains fixed. 
Now, when we continue to $t=1$, the points of the fiber continue 
to move counterclockwise around $y=0$. When we reach $t=1$ and 
reach back the fiber $x=1$, the points in the fiber made a $180^{\circ}$ rotation 
counterclockwise around $y=0 \in \C$ from their corresponding points in 
the initial fiber $x=1$. Hence, we have: 

\begin{thm}
The action of local braid monodromy 
of the point presented locally by the equation: $y(y^2+x)(y^2-x)=0$ is a 
$180^{\circ}$ rotation counterclockwise of the four points 
(correspond to $y=1,-1,i,-i$) around the point corresponds to $y=0$, 
as shown schematically in Figure \ref{bm_type_h}.
 
\begin{figure}[h]
\epsfysize=2.5cm  
\centerline{\epsfbox{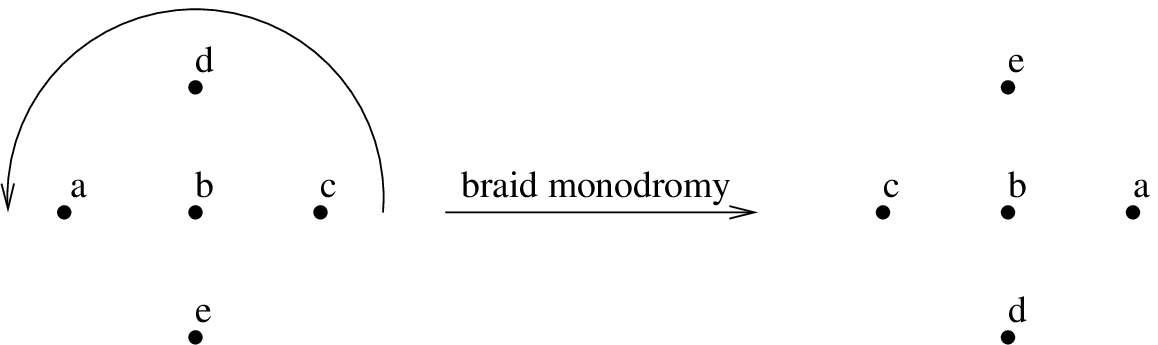}}  
\caption{The local braid monodromy of the singularity $y(y^2+x)(y^2-x) = 0$}\label{bm_type_h}  
\end{figure}  
\end{thm}
 
\begin{figure}[h]
\epsfysize=12cm  
\centerline{\epsfbox{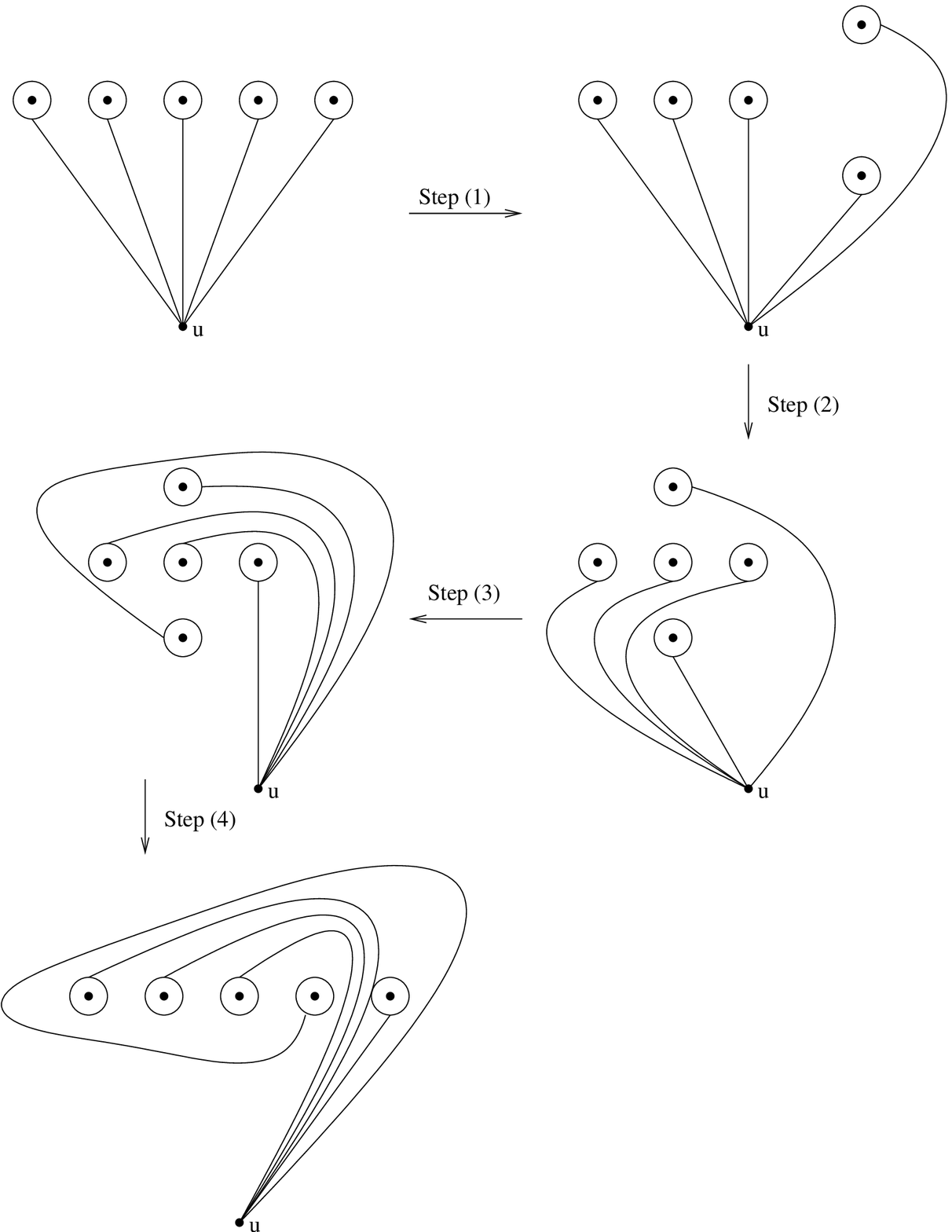}}  
\caption{The g-base obtained from the standard g-base by the action of the local braid monodromy}\label{gbase_type_h}  
\end{figure}

Now, we want to compute the induced relations of this singular point. By van 
Kampen's Theorem, one should compute the g-base obtained by applying the action induced
by the local braid monodromy of the singular point on the standard g-base. Since we have 
two complex points in the fiber  before the action of the braid monodromy and after it,
we have to start by rotating the two rightmost points by $90^{\circ}$ 
counterclockwise, for representing the two complex points (see Step (1) in 
Figure \ref{gbase_type_h}). Then, we move the two complex points to be $i$ and 
$-i$ (see Step (2) there). 
Now, we apply the action of the local braid monodromy (Step (3)),
and then we return the two new complex points to the right side, and return them to the 
real axis by rotating them by $90^{\circ}$ clockwise (Step (4)). 

\begin{rem}\label{rem_clockwise}
One should notice here that when we return the points back to the real axis in Step (4), 
we rotate the points clockwise. The reason 
is that we want to compute the pure action of the braid on the g-base, and 
the clockwise rotation cancels the counterclockwise rotation which we have performed 
in Step (1).
In the global computations, we indeed rotate always counterclockwise, since we work 
there in a model (for simplifying the computations), and we can define the model 
as we want. We will get the same results in the global computations even if we rotate 
clockwise the complex points into real points.    
\end{rem}

Now, by van Kampen's Theorem, we get the following induced relations from 
the new g-base (where $\{ x_1,x_2,x_3,x_4,x_5 \}$ are the generators of the standard 
g-base):
\begin{enumerate}
\item $x_1 = x_4 x_3 x_4^{-1}$ 
\item $x_2 = x_4 x_3 x_2 x_3^{-1} x_4^{-1}$
\item $x_3 = x_4 x_3 x_2 x_1 x_2^{-1} x_3^{-1} x_4^{-1}$
\item $x_4 = x_5$
\item $x_5 = x_5 x_4 x_3 x_2 x_1 x_4 x_1^{-1} x_2^{-1} x_3^{-1} x_4^{-1} x_5^{-1}$  
\end{enumerate}

Relation (2) is equivalent to $x_4 x_3 x_2 = x_2 x_4 x_3$.
Substituting $x_1$ by $x_4 x_3 x_4^{-1}$ in Relation (3) yields $x_3 = x_4 x_3 x_2 x_4 x_3 x_4^{-1} x_2^{-1} x_3^{-1} x_4^{-1}$, which is equivalent to $x_4 x_3 x_2 x_4 x_3 = x_3 x_4 x_3 x_2 x_4$. By Relation (2), it can also be written as $x_4 x_3 x_2 x_4 x_3 = x_3 x_2 x_4 x_3 x_4$.

We will show that Relation (5) is redundant. First, we can cancel $x_5$, so we get:
$x_4 x_3 x_2 x_1 x_4 = x_5 x_4 x_3 x_2 x_1$. Since $x_4=x_5$, 
we can cancel another $x_5$,
and by substituting $x_1$ by $x_4 x_3 x_4^{-1}$ and some simplifications we get 
$x_3 x_2 x_4 x_3 x_4 = x_4 x_3 x_2 x_4 x_3$, which is equal to Relation (3). 
Hence, Relation (5) is redundant.

Therefore, we have the following corollary:

\begin{cor}
The induced relations for the singular point presented locally by the equation 
$y(y^2+x)(y^2-x)=0$ are:
\begin{enumerate}
\item $x_4 x_3 x_2 = x_2 x_4 x_3$
\item $x_3 x_2 x_4 x_3 x_4 = x_4 x_3 x_2 x_4 x_3$
\item $x_1 = x_4 x_3 x_4^{-1}$
\item $x_4 = x_5$
\end{enumerate}
where $\{ x_1,x_2,x_3,x_4,x_5 \}$ are the generators of the standard g-base.
\end{cor}

Now, we compute the corresponding Lefschetz diffeomorphism: 
\begin{cor}
The Lefschetz diffeomorphism of the singular point presented locally by 
$y(y^2+x)(y^2-x)=0$ is a $90^{\circ}$ rotation counterclockwise of the
four points around the fixed point $y=0 \in \C$. 
\end{cor}

The reason for this corollary is since the Lefschetz diffeomorphism is obtained
by computing the action along the loop mentioned above from $t=\frac{1}{2}$ 
to $t=1$. Along this interval, the four points (except for $y=0$) made a 
$90^{\circ}$ rotation counterclockwise.

\subsection{Fourth type}\label{3comps_type4}
 
This type is slightly different from the previous three types. In the previous types,
the line intersects the tangency point, but was not tangented to it. 
In this type, the additional line and the tangency point have a common tangent.

\begin{figure}[h]
\epsfysize=2.5cm  
\centerline{\epsfbox{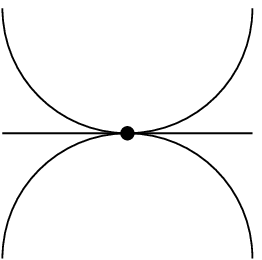}}  
\caption{The singularity of $y(y+x^2)(y-x^2) = 0$ at $(0,0)$}\label{type_d}  
\end{figure}

The local equation of the singularity of this type is 
$y(y+x^2)(y-x^2) = 0$  (see Figure \ref{type_d}). 

We take a loop around $x=0$ and look what happens
to the points of the curve $y(y+x^2)(y-x^2) = 0$ in the fibers 
over this loop. 

Let $x=e^{2\pi i t}$ where $0 \leq t \leq 1$. For $t=0$, we have that $x=1$ and 
the points of the curve in the fiber over $x=1$ are $y=-1,0$ and $y=1$. 
For $t=\frac{1}{2}$, we have $x=-1$, and we get again that the points of the curve 
in the fiber over $x=-1$ are $y=-1,0$ and $y=1$.  One can see that the point $y=0$ is 
fixed while $t$ is changing from $0$ to $\frac{1}{2}$, and the other two points,
$y=1$ and $y=-1$, do a counterclockwise full-twist around the point $y=0$.
Hence, when we continue to $t=1$, we get that the point $y=1$ and $y=-1$ do two 
counterclockwise full-twists around the fixed point $y=0$.

\begin{figure}[h]
\epsfysize=4cm  
\centerline{\epsfbox{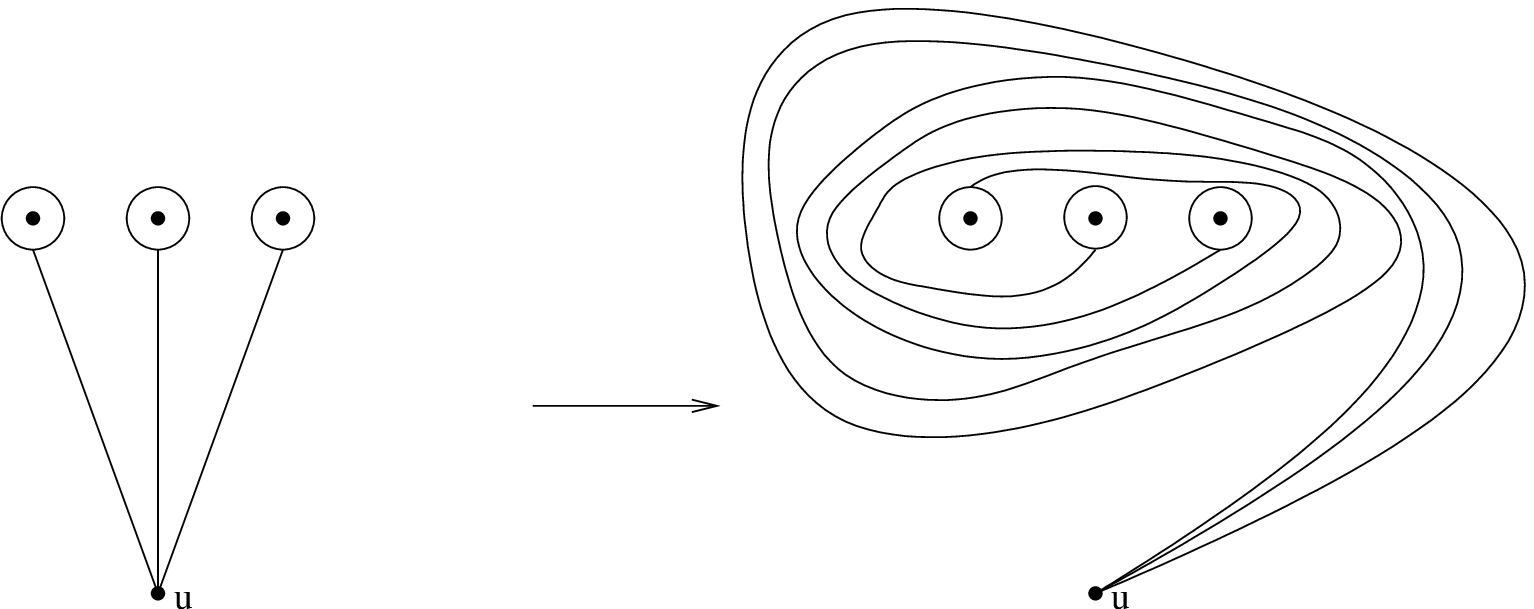}}  
\caption{The g-base obtained from the standard g-base by the action of the local braid monodromy}\label{gbase_type_d}  
\end{figure}

As before, by van Kampen's Theorem, one should compute the g-base 
obtained by applying the action induced
by the local braid monodromy of the singular point on the standard g-base. 
This is shown in Figure \ref{gbase_type_d}.

Now, by van Kampen's Theorem, we get the following induced relations from 
the new g-base (where $\{ x_1,x_2,x_3 \}$ are the generators of the standard g-base):
\begin{enumerate}
\item $x_1 = x_3 x_2 x_1 x_3 x_2 x_1 x_2^{-1} x_3^{-1} x_1^{-1} x_2^{-1} x_3^{-1}$ 
\item $x_2 = x_3 x_2 x_1 x_3 x_2 x_1 x_2 x_1^{-1} x_2^{-1} x_3^{-1} x_1^{-1} x_2^{-1} x_3^{-1}$
\item $x_3 = x_3 x_2 x_1 x_3 x_2 x_1 x_3 x_1^{-1} x_2^{-1} x_3^{-1} x_1^{-1} x_2^{-1} x_3^{-1}$
\end{enumerate}

From Relation (1), one gets $(x_3 x_2 x_1)^2 = (x_1 x_3 x_2)^2$. Relation (2) can be written
$x_2 (x_3 x_2 x_1)^2 = (x_3 x_2 x_1)^2 x_2$. By Relation (1), it becomes
$x_2 (x_1 x_3 x_2)^2 = (x_1 x_3 x_2)^2 x_2$. By cancelling  the rightmost $x_2$, we get:
$(x_2 x_1 x_3)^2 = (x_1 x_3 x_2)^2$.

Relation (3) is redundant, since by cancelling $x_3$, one gets:
$(x_3 x_2 x_1)^2 = (x_2 x_1 x_3)^2$, which is already known.

To summarize, we have proved the following result:
\begin{thm}
The local braid monodromy of the singularity presented locally by the equation 
$y(y+x^2)(y-x^2)=0$ is: the three points do two 
counterclockwise generalized full-twists.

The induced relations of this point are: 
$$(x_3 x_2 x_1)^2 = (x_1 x_3 x_2)^2 = (x_2 x_1 x_3)^2,$$
where $\{ x_1,x_2,x_3 \}$ are the generators of the standard g-base.
\end{thm}

The corresponding Lefschetz diffeomorphism is presented in the following corollary: 
\begin{cor}
The Lefschetz diffeomorphism of the singular point presented locally by 
$y(y+x^2)(y-x^2)=0$ is a full-twist of all the three points. 
\end{cor}

The reason for that corollary is since the Lefschetz diffeomorphism is obtained
by going along the loop mentioned above from $t=\frac{1}{2}$ 
to $t=1$. In this interval, the three points perform a full-twist.

\bigskip

The arguments of this section can be generalized easily to an arbitrary number 
of tangented components at the same point with a common tangent:

\begin{cor}
Let $P$ be a singular point which consists of $n$ smooth branches 
with a common tangent at $P$. Let $x_i$ be the generator which corresponds
to the $i$th branch. Then, the local braid monodromy is two generalized full-twists 
of all the segment from $1$ to $n$. Moreover, the relations induced by 
this singular point are:
$$(x_n x_{n-1} \cdots x_1)^2 = (x_{n-1} \cdots x_1 x_n)^2 = \cdots = (x_1 x_n \cdots x_2)^2.$$
\end{cor}

Notice that if the tangented components have a higher order of tangency (i.e. 
the higher derivatives of the components are also equal up to the $k$th derivative,
for a given $k$), then the braid monodromy will consist of an higher power of 
a full-twist, and a higher exponent will appear in the relations.   

\section{Singular points with four components}\label{4_comps}

In this section, we will compute the local braid monodromy of five cases 
of singular points consist of four components, where at most two 
of them are conics. 
Then, we will compute the relations induced from these singular 
points by the classical van Kampen Theorem
(see \cite{vK} and Section \ref{vK-thm} here). 

A singular point consists of four components in a conic-line arrangements with 
two tangented conics and up to two additional lines can be of the following 
two types: take the tangency point of the two conics, and add two lines in 
the following two ways: 
one way is to add one line which will be tangented to both conics, 
and the second line will intersect both conics (and the line) at the 
singular point. The second way will be to add two intersecting lines 
at the tangency point of the two conics.

\subsection{Two tangented conics with a tangent line and an intersecting line}

As before, there are three possibilities for the intersecting line: it can 
locally presented as the line $y=2x$, $y=-2x$ or $x=0$. In spite of the
fact that one can move from the first type to the second type by rotating 
the line, so locally the singularities are equivalent, 
but from the global point of view 
these singularities induced different relations in the global fundamental group.
     
\subsubsection{First type}\label{4comps_type1}

\begin{figure}[h]
\epsfysize=2.5cm  
\centerline{\epsfbox{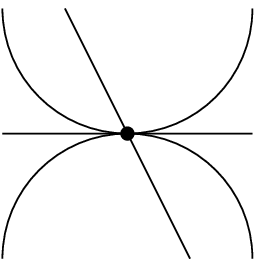}}  
\caption{The singularity of $y(2x+y)(y+x^2)(y-x^2) = 0$ at $(0,0)$}\label{type_e}  
\end{figure}

The local equation of the singularity of the first type is 
$y(2x+y)(y+x^2)(y-x^2) = 0$ (see Figure \ref{type_e}). 

For computing the braid monodromy, we take a loop around $x=0$ and 
look what happens to the points of the curve $y(2x+y)(y+x^2)(y-x^2) = 0$ 
in the fibers above this loop. 
 
Let $x=e^{2\pi i t}$ where $0 \leq t \leq 1$. For $t=0$ we have that 
$x=1$ and the points of the curve in the fiber over $x=1$ are $y=1,0,-1$ 
and $y=-2$. For $t=\frac{1}{2}$, we have $x=-1$, 
and the points of the curve in the fiber over $x=-1$ are $y=1,0,-1$ and $y=2$.  
The point $y=1$
in the fiber $x=-1$ corresponds to the point $y=1$ in the fiber $x=1$, and similarly
the point $y=-1$ in the fiber $x=-1$ corresponds to the point $y=-1$ 
in the fiber $x=1$. The point $y=0$ is fixed, and the point $y=2$ 
in the fiber $x=-1$ corresponds to the point $y=-2$ in the fiber $x=1$. 
Hence, we get that from $t=0$ to $t=\frac{1}{2}$, the points $y=-1,0$
and $y=1$ do a counterclockwise generalized full-twist 
and the point $y=-2$ does 
a counterclockwise half-twist around these three points. 
If we continue to $t=1$, the three points $y=-1,0,1$ perform 
two counterclockwise generalized full-twists and the point $y=-2$ does 
a counterclockwise full-twist around the other three points. 

\begin{figure}[h]
\epsfysize=8cm  
\centerline{\epsfbox{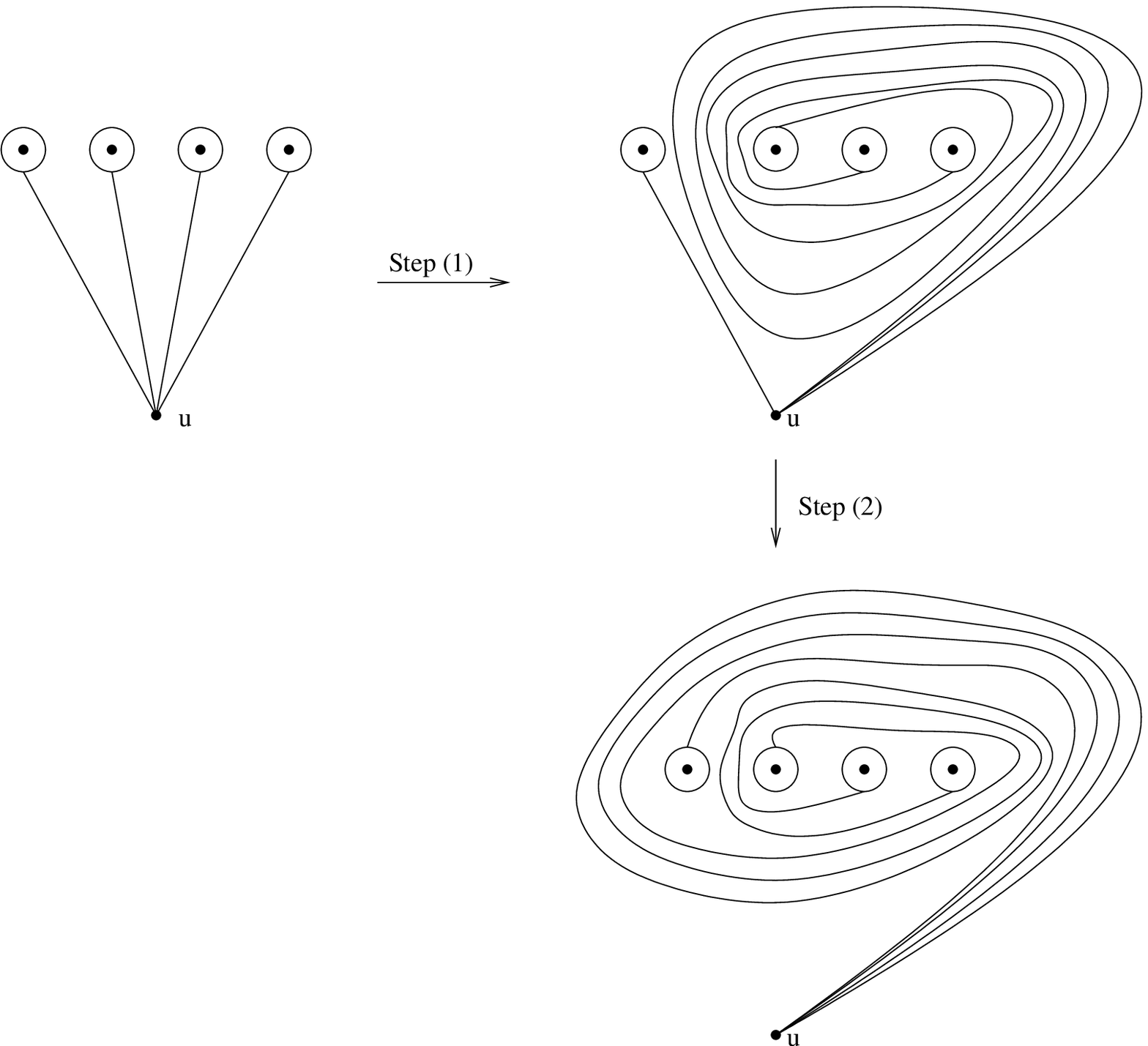}}  
\caption{The g-base obtained from the standard g-base by the action of the local braid monodromy}\label{gbase_type_e}  
\end{figure}

For computing the induced relations, we use the van Kampen Theorem. 
Here, one should compute the g-base obtained by applying the action induced
by the local braid monodromy of the singular point on the standard g-base. 
This is shown in Figure \ref{gbase_type_e}: in the first step, we perform two 
counterclockwise generalized full-twists of points $2,3$ and $4$ which represent 
$y=-1,0$ and $y=1$ respectively. In the second step, 
we perform a counterclockwise full-twist of point $1$ (which represents 
the point $y=-2$) around the other three points.

By van Kampen's Theorem, we get the following induced relations
(where $\{ x_1,x_2,x_3,x_4 \}$ are the generators of the standard g-base):
\begin{enumerate}
\item $x_1 = x_4 x_3 x_2 x_1 x_2^{-1} x_3^{-1} x_4^{-1}$  
\item $x_2 = x_4 x_3 x_2 x_1 x_4 x_3 x_2 x_3^{-1} x_4^{-1} x_1^{-1} x_2^{-1} x_3^{-1} x_4^{-1}$
\item $x_3 = x_4 x_3 x_2 x_1 x_4 x_3 x_2 x_3 x_2^{-1} x_3^{-1} x_4^{-1} x_1^{-1} x_2^{-1} x_3^{-1} x_4^{-1}$
\item $x_4 = x_4 x_3 x_2 x_1 x_4 x_3 x_2 x_4 x_2^{-1} x_3^{-1} x_4^{-1} x_1^{-1} x_2^{-1} x_3^{-1} x_4^{-1}$
\end{enumerate}

From Relation (1), we get that $x_1 x_4 x_3 x_2 = x_4 x_3 x_2 x_1$. 
Relations (2) and (3) become 
$x_4 x_3 x_2 x_1 x_4 x_3 x_2 = x_2 x_4 x_3 x_2 x_1 x_4 x_3$ and 
$x_4 x_3 x_2 x_1 x_4 x_3 x_2 x_3 = x_3 x_4 x_3 x_2 x_1 x_4 x_3 x_2$ 
respectively. By Relation (2), Relation (3) can be written as:
$$x_4 x_3 x_2 x_1 x_4 x_3 x_2 x_3 = x_3 x_2 x_4 x_3 x_2 x_1 x_4 x_3,$$ and hence
we get from Relations (2) and (3):
$$x_4 x_3 x_2 x_1 x_4 x_3 x_2 = x_3 x_2 x_4 x_3 x_2 x_1 x_4 = x_2 x_4 x_3 x_2 x_1 x_4 x_3.$$ 

We will show now  that Relation (4) is redundant. 
From Relation (4), one can get:
$$x_3 x_2 x_1 x_4 x_3 x_2 x_4 = x_4 x_3 x_2 x_1 x_4 x_3 x_2.$$
Using Relation (1), this relation is equal to Relation (3), and hence 
it is redundant. 

Hence, we have:

\begin{thm}
The local braid monodromy of the singularity presented locally by the equation 
$y(2x+y)(y+x^2)(y-x^2)=0$ is: three points (correspond to $y=1,0$ and $y=-1$) do two 
counterclockwise generalized full-twists and the fourth point (corresponds to $y=-2$) do 
a counterclockwise full-twist around them. 

The induced relations from this singular point are:
$$x_1 x_4 x_3 x_2 = x_4 x_3 x_2 x_1 \quad ; \quad (x_4 x_3 x_2)^2 x_1 = x_3 x_2 x_4 x_3 x_2 x_1 x_4 = x_2 x_4 x_3 x_2 x_1 x_4 x_3,$$ 
where $\{ x_1,x_2,x_3,x_4 \}$ are the generators of the standard g-base. 
\end{thm}

One should notice here that if we delete the generator $x_1$ which corresponds 
to the line that intersects the three tangented components, we get the following 
set of relations:
$$(x_4 x_3 x_2)^2 = (x_3 x_2 x_4)^2 = (x_2 x_4 x_3)^2,$$
as expected (see Section \ref{3comps_type4}).
Similarly, if we delete the generator $x_3$ which corresponds to the tangent line, 
we get the following set of relations:
$$x_1 x_4 x_2 = x_4 x_2 x_1 \quad ; \quad x_4 x_2 x_4 x_2 x_1 = x_2 x_4 x_2 x_1 x_4,$$ 
again as expected (see Section \ref{3comps_type1}).

\medskip

Since the Lefschetz diffeomorphism of the singular point 
is obtained by computing
the action only on half of the unit circle (from $t=\frac{1}{2}$ to $t=1$), 
we have the following corollary:

\begin{cor}
The Lefschetz diffeomorphism of the singular point presented locally by
$y(2x+y)(y+x^2)(y-x^2) = 0$ is a counterclockwise generalized full-twist 
of the points correspond to $y=1,0$ and $y=-1$, and a half-twist 
of the point corresponds to $y=2$ (which becomes $y=-2$).
\end{cor}

\subsubsection{Second type}\label{4comps_type2}
\begin{figure}[h]
\epsfysize=2.5cm  
\centerline{\epsfbox{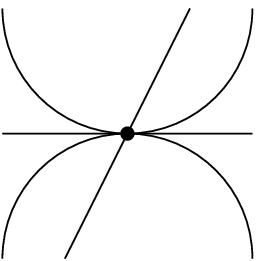}}  
\caption{The singularity of $y(2x-y)(y+x^2)(y-x^2) = 0$ at $(0,0)$}\label{type_f}  
\end{figure}

The local equation of the singularity of the second type is 
$y(2x-y)(y+x^2)(y-x^2) = 0$  (see Figure \ref{type_f}). 

Applying almost the same computations we have performed for the previous type, 
we get:

\begin{thm}
The local braid monodromy of the singularity presented locally by the equation 
$y(2x-y)(y+x^2)(y-x^2)=0$ is: the points $y=1,0,-1$ do two 
counterclockwise generalized full-twists, and the point $y=2$ does a 
counterclockwise full-twist around the points $y=1,0,-1$ together.

The corresponding induced relations are:    
$$x_4 x_3 x_2 x_1 = x_3 x_2 x_1 x_4 \quad ; \quad (x_3 x_2 x_1)^2 x_4 = x_2 x_1 x_3 x_2 x_1 x_4 x_3 = x_1 x_3 x_2 x_1 x_4 x_3 x_2,$$ 
where $\{ x_1,x_2,x_3,x_4 \}$ are the generators of the standard g-base.
\end{thm}

The Lefschetz diffeomorphism is presented in the following corollary.

\begin{cor}
The Lefschetz diffeomorphism of the singular point presented locally by
$y(2x-y)(y+x^2)(y-x^2) = 0$ is a counterclockwise generalized full-twist 
of the points correspond to $y=1,0$ and $y=-1$, and a half-twist 
of the point corresponds to $y=-2$ (which becomes $y=2$).
\end{cor}

\subsubsection{Third type}\label{4comps_type3}
\begin{figure}[h]
\epsfysize=2.5cm  
\centerline{\epsfbox{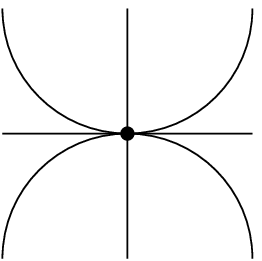}}  
\caption{The singularity of $xy(y+x^2)(y-x^2) = 0$ at $(0,0)$}\label{type_j}  
\end{figure}

The local equation of the singularity of the third type is 
$xy(y+x^2)(y-x^2) = 0$  (see Figure \ref{type_j}). 

For computing the braid monodromy in this case, we will use the following trick:
the three components in the middle $x(y+x^2)(y-x^2)=0$ can be thought 
for a moment as a ``thick'' line which is perpendicular to the other line.
After this observation, we have two intersecting lines, 
which their local braid monodromy is a counterclockwise full-twist (even though 
we have here a thick line $x=0$, which is vertical, one can rotate it a bit for 
getting two ``usual'' intersecting lines, without changing the braid monodromy).

The local braid monodromy induced by the ``thick'' line itself (which consists 
of the curve $x(y+x^2)(y-x^2)=0$) has already been computed (Section 
\ref{3comps_type3}): 
we got there that the local braid monodromy is a $180^{\circ}$ counterclockwise 
rotation of four points around one fixed point at the origin. 

Hence, the action of the local braid monodromy of the whole singular point
can be summarized as follows: 

\begin{thm}
The local braid monodromy of the singular point presented locally by the 
equation $xy(y+x^2)(y-x^2)=0$ is: first perform a $180^{\circ}$ counterclockwise 
rotation of the points correspond to $y=1,-1,i,-i$ around $y=0$, 
and then the point corresponds to $y=2$ performs a 
full-twist with a block which consists of all the other points (see Figure 
\ref{bm_type_j}).

\begin{figure}[h]
\epsfysize=2.8cm  
\centerline{\epsfbox{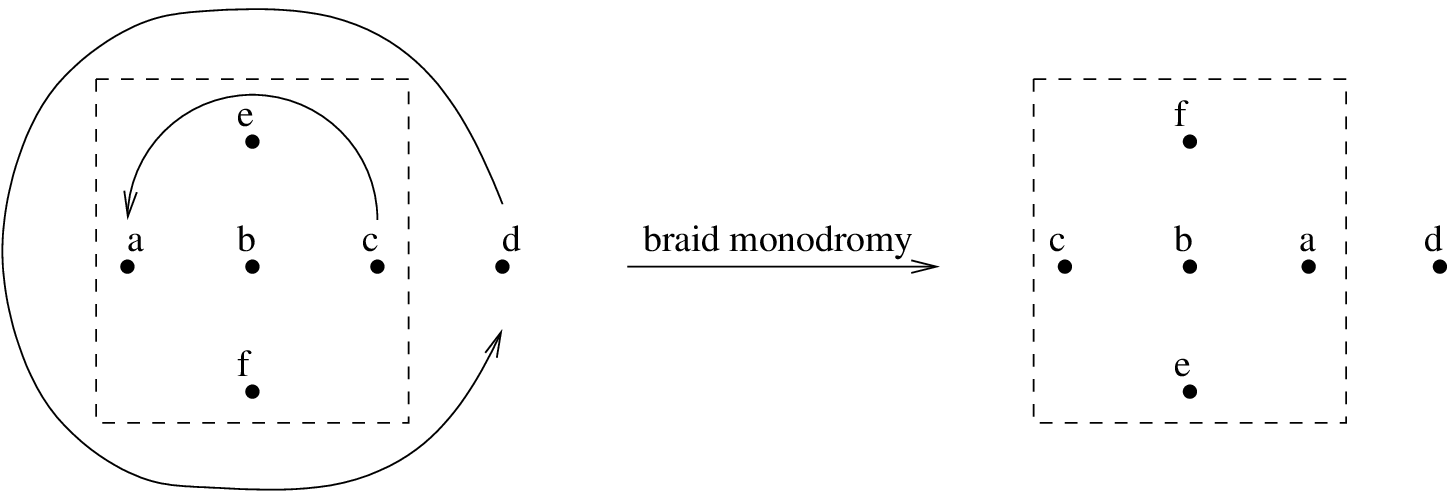}}  
\caption{The action of the local braid monodromy}\label{bm_type_j}  
\end{figure}  
\end{thm}

\begin{figure}[h]
\epsfysize=12cm  
\centerline{\epsfbox{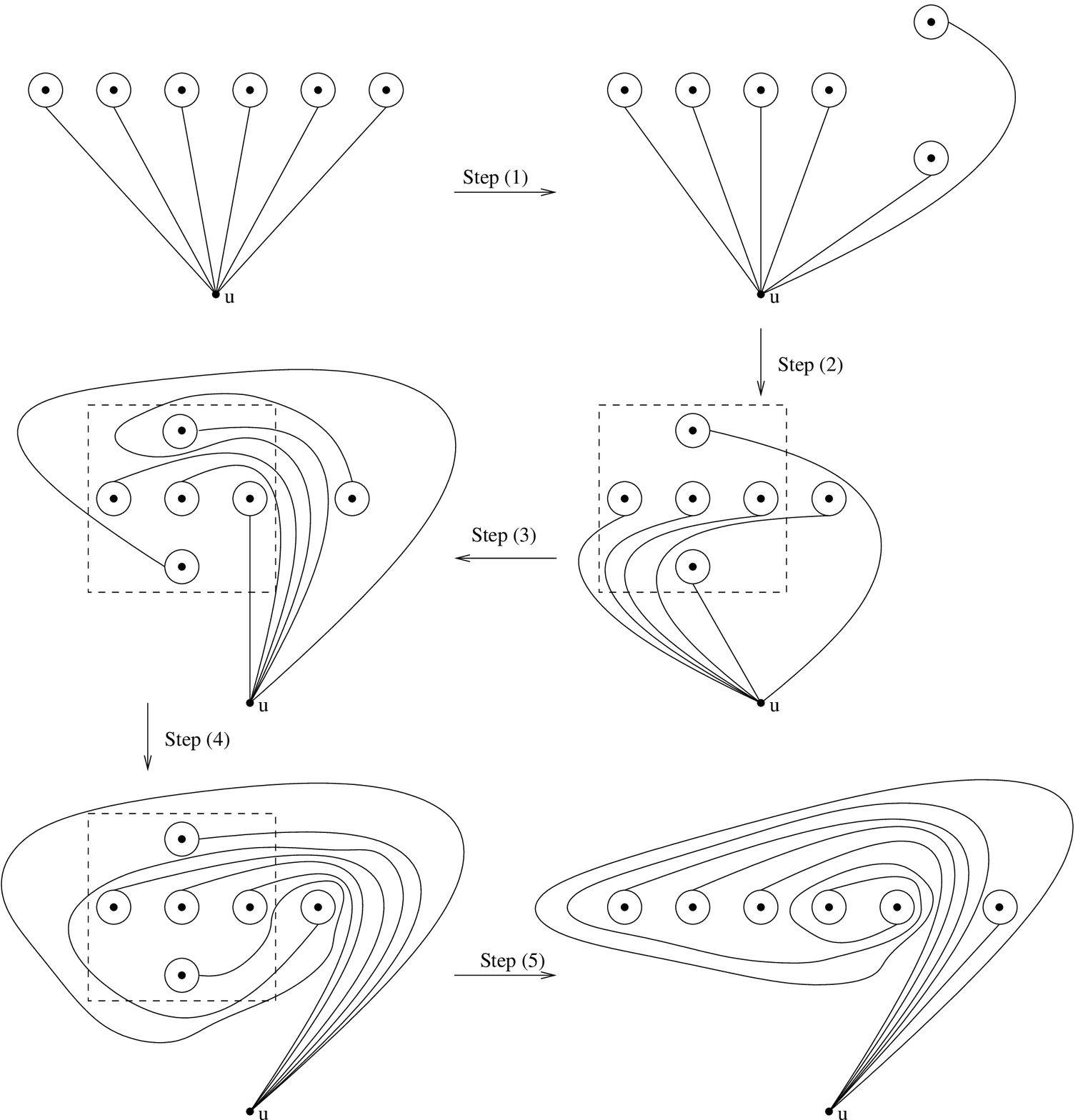}}  
\caption{The g-base obtained from the standard g-base by the action of the local braid monodromy}\label{gbase_type_j}  
\end{figure}

Now, we want to compute the induced relations of this singular point. 
According to van Kampen's Theorem, one should compute the g-base 
obtained by applying the action induced by the local braid monodromy 
of the singular point on the standard g-base. Since we have 
two complex points in the fiber  before the action 
of the braid monodromy and after it, we have to start by rotating 
the two rightmost points by $90^{\circ}$ 
counterclockwise, for representing the two complex points (see Step (1) in 
Figure \ref{gbase_type_j}). Then, we move the two complex points to be over 
the second to the left real point (see Step (2) there). 
Now, we apply the action of the local braid monodromy (Steps (3) and (4)),
and then we return the two new complex points to the right side, 
and return them to the real axis by rotating them clockwise by 
$90^{\circ}$ (Step (5), see Remark \ref{rem_clockwise}). 

Now, by van Kampen's Theorem, we get the following induced relations from 
the new g-base (where $\{ x_1,x_2,x_3,x_4,x_5,x_6 \}$ are the generators of the standard g-base):
\begin{enumerate}
\item $x_1=x_5 x_4 x_3 x_4^{-1} x_5^{-1}$
\item $x_2=x_5 x_4 x_3 x_2 x_3^{-1} x_4^{-1} x_5^{-1}$
\item $x_3=x_5 x_4 x_3 x_2 x_1 x_2^{-1} x_3^{-1} x_4^{-1} x_5^{-1}$
\item $x_4=x_5 x_4 x_3 x_2 x_1 x_5 x_4 x_5^{-1} x_1^{-1} x_2^{-1} x_3^{-1} x_4^{-1} x_5^{-1}$
\item $x_5=x_6$
\item $x_6=x_6 x_5 x_4 x_3 x_2 x_1 x_5 x_4 x_5 x_4^{-1} x_5^{-1} x_1^{-1} x_2^{-1} x_3^{-1} x_4^{-1} x_5^{-1} x_6^{-1}$ 
\end{enumerate}

Relation (2) can be written as: $x_5 x_4 x_3 x_2 = x_2 x_5 x_4 x_3$.
Using Relation (1), Relation (3) becomes: 
$x_3 x_5 x_4 x_3 x_2 x_5 x_4 = x_5 x_4 x_3 x_2 x_5 x_4 x_3$.

By Relation (1) and some cancellations, Relation (4) becomes: 
$$x_4=x_5 x_4 x_3 x_2 x_5 x_4 x_3 x_4 x_3^{-1} x_4^{-1} x_5^{-1} x_2^{-1} x_3^{-1} x_4^{-1} x_5^{-1}.$$
By Relation (3), we get:
$$x_4=x_3 x_5 x_4 x_3 x_2 x_5 x_4 x_5^{-1} x_2^{-1} x_3^{-1} x_4^{-1} x_5^{-1} x_3^{-1},$$
and hence we have: $x_4 x_3 x_5 x_4 x_3 x_2 x_5 = x_3 x_5 x_4 x_3 x_2 x_5 x_4$.

Now, we will show that Relation (6) is redundant. First, this relation 
can be written as: $x_5 x_4 x_3 x_2 x_1 x_5 x_4 = x_4 x_3 x_2 x_1 x_5 x_4 x_5$. 
By Relation (1) and some cancellations, we have:
$x_5 x_4 x_3 x_2 x_5 x_4 x_3 = x_4 x_3 x_2 x_5 x_4 x_3 x_5$.
By Relation (3), we have: 
$x_3 x_5 x_4 x_3 x_2 x_5 x_4 = x_4 x_3 x_2 x_5 x_4 x_3 x_5$.
By Relation (2), this relation is equal to Relation (4), and hence Relation (6) is 
redundant.

Therefore, we get the following result for the set of relations for the singular point:

\begin{cor}
The singular point presented locally by $xy(y+x^2)(y-x^2) = 0$ 
has the following set of induced relations: 
\begin{enumerate}
\item $x_5 x_4 x_3 x_2 = x_2 x_5 x_4 x_3$
\item $x_3 x_5 x_4 x_3 x_2 x_5 x_4 = x_5 x_4 x_3 x_2 x_5 x_4 x_3 = x_4 x_3 x_5 x_4 x_3 x_2 x_5$
\item $x_1 = x_5 x_4 x_3 x_4^{-1} x_5^{-1}$
\item $x_5 = x_6$
\end{enumerate}
where $\{ x_1,x_2,x_3,x_4,x_5,x_6 \}$ are the generators of the standard g-base.
\end{cor}

One should notice here that if we delete the generator $x_2$ which corresponds 
to the vertical line that intersects the three tangented components, 
we get the following set of relations:
$$(x_5 x_4 x_3)^2 = (x_4 x_3 x_5)^2 = (x_3 x_5 x_4)^2,$$
as expected (see Section \ref{3comps_type4}).
Similarly, if we delete the generator $x_4$ which corresponds to the tangent line, 
we get the relations of the singular point in Section \ref{3comps_type3}
as expected.

\medskip

Since the Lefschetz diffeomorphism of the singular point 
is obtained by computing
the action only on half of the unit circle (from $t=\frac{1}{2}$ to $t=1$), 
we have the following corollary:

\begin{cor}
The Lefschetz diffeomorphism of the singular point presented locally by
$xy(y+x^2)(y-x^2) = 0$ is a $90^{\circ}$ counterclockwise rotation of the 
points correspond to $y=1,-1,i,-i$ around $y=0$, 
and then the point corresponds to $y=2$ do a 
half-twist with a block which consists of all the other points.
\end{cor}

\subsection{Two tangented conics with two intersecting lines}

Although the following two types are almost a rotation of each other, we need them both, 
since the second type includes also ``hidden'' branch points inside the 
singularity, as we had in Section \ref{3comps_type3}.

\subsubsection{First type}\label{4comps_type4}
\begin{figure}[h]
\epsfysize=2.5cm  
\centerline{\epsfbox{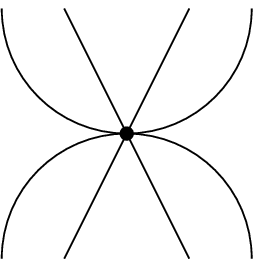}}  
\caption{The singularity of $(2x+y)(2x-y)(y+x^2)(y-x^2) = 0$ at $(0,0)$}\label{type_k}  
\end{figure}

The local equation of the singularity of the first type is 
$(2x+y)(2x-y)(y+x^2)(y-x^2) = 0$ (see Figure \ref{type_k}). 

For computing the braid monodromy, we take a loop around $x=0$ and look what happens
to the points of the curve $(2x+y)(2x-y)(y+x^2)(y-x^2) = 0$ in the fibers 
over this loop. 
 
Let $x=e^{2\pi i t}$ where $0 \leq t \leq 1$. For $t=0$ we have that $x=1$ and 
the points of the curve in the fiber over $x=1$ are $y=1,-1,2$ and $y=-2$. 
For $t=\frac{1}{2}$, we have $x=-1$, and the points of the curve 
in the fiber over $x=-1$ are $y=1,-1,-2$ and $y=2$.  The point $y=1$
in the fiber $x=-1$ corresponds to the point $y=1$ in the fiber $x=1$, and similarly
the point $y=-1$ in the fiber $x=-1$ corresponds to the point $y=-1$ 
in the fiber $x=1$. The point $y=2$ in the fiber $x=-1$ corresponds to 
the point $y=-2$ in the fiber $x=1$, and similarly the point $y=-2$ 
in the fiber $x=-1$ corresponds to the point $y=2$ in the fiber $x=1$. 
Hence, we get that from $t=0$ to $t=\frac{1}{2}$, the points $y=1$
and $y=-1$ do a counterclockwise full-twist and the points $y=2$ and $y=-2$ do 
a counterclockwise half-twist around the points $y=1,-1$. 
If we continue to $t=1$, the points $y=1$
and $y=-1$ do two counterclockwise full-twists and the points $y=2$ and $y=-2$ do 
a counterclockwise full-twist around the points $y=1,-1$ together. 

Hence we have the following result:
\begin{thm}
The local braid monodromy of the singularity presented locally by the equation 
$(2x+y)(2x-y)(y+x^2)(y-x^2) = 0$ is: two points (correspond to $y=1$ and $y=-1$) do two 
counterclockwise full-twists and the two other points (corresponds to $y=2$ and $y=-2$) do 
a counterclockwise full-twist around them. 
\end{thm}

\begin{figure}[h]
\epsfysize=8cm  
\centerline{\epsfbox{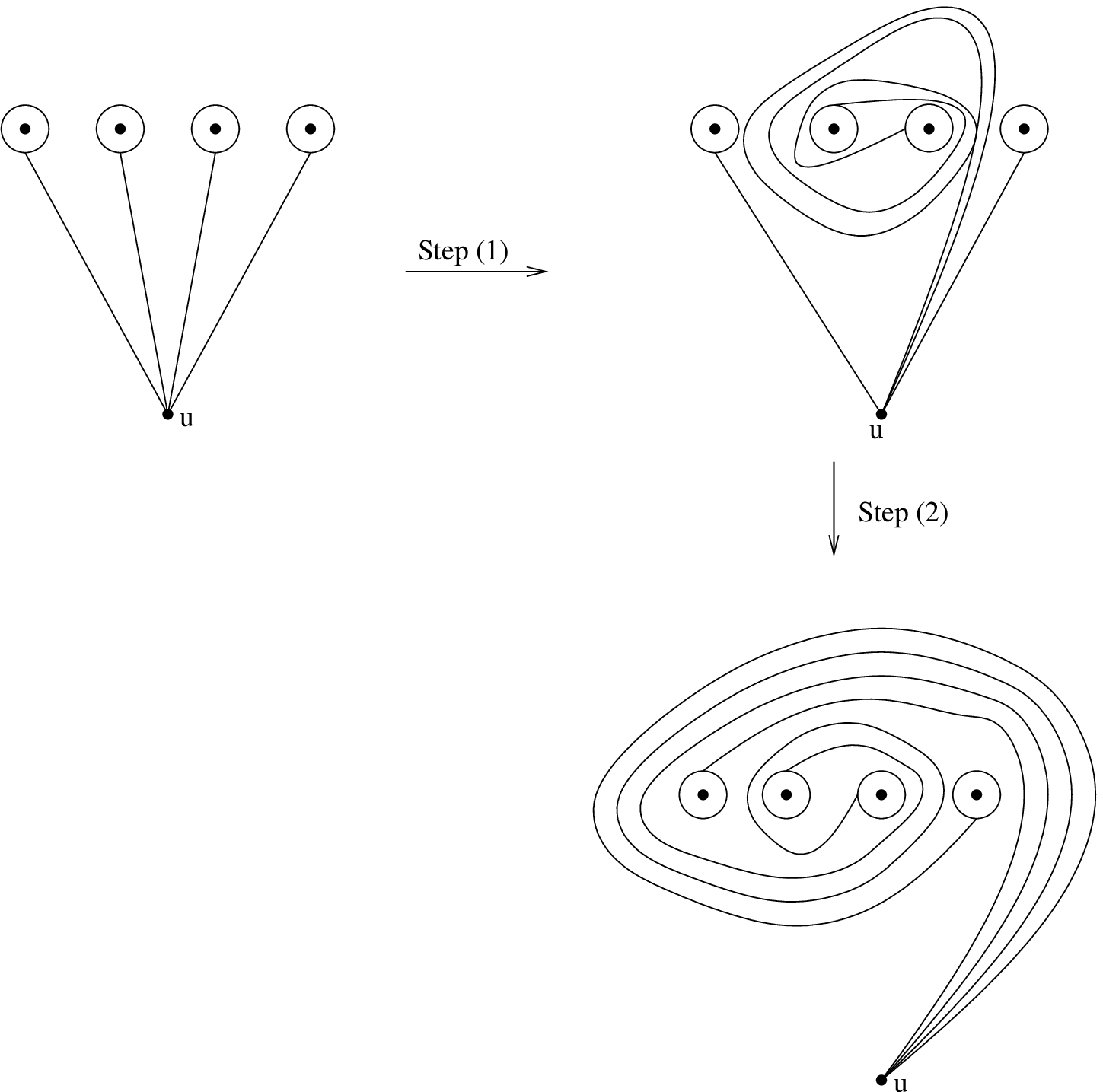}}  
\caption{The g-base obtained from the standard g-base by the action of the local braid monodromy}\label{gbase_type_k}  
\end{figure}

Now, we want to compute the induced relations from the braid monodromy 
of the singular point. By van Kampen's Theorem, we compute the g-base obtained by applying the action induced
by the local braid monodromy of the singular point on the standard g-base. 
This is shown in Figure \ref{gbase_type_k}: in the first step, we perform two 
counterclockwise full-twists of point $2$ and $3$ which represent 
$y=-1$ and $y=1$ respectively. In the second step, 
we perform a counterclockwise full-twist of points $1$ and $4$ (which represent 
$y=2$ and $y=-2$) around the points $2$ and $3$.

Now, by van Kampen's Theorem, we get the following induced relations from 
the new g-base (where $\{ x_1,x_2,x_3,x_4 \}$ are the generators 
of the standard g-base):

\begin{enumerate}
\item $x_1 = x_4 x_3 x_2 x_1 x_2^{-1} x_3^{-1} x_4^{-1}$ 
\item $x_2 = x_4 x_3 x_2 x_1 x_3 x_2 x_3^{-1} x_1^{-1} x_2^{-1} x_3^{-1} x_4^{-1}$
\item $x_3 = x_4 x_3 x_2 x_1 x_3 x_2 x_3 x_2^{-1} x_3^{-1} x_1^{-1} x_2^{-1} x_3^{-1} x_4^{-1}$
\item $x_4 = x_4 x_3 x_2 x_1 x_4 x_1^{-1} x_2^{-1} x_3^{-1} x_4^{-1}$
\end{enumerate}

From Relation (1), we get $x_4 x_3 x_2 x_1 = x_1 x_4 x_3 x_2$. From Relation (4), we get
$x_4 x_3 x_2 x_1 = x_3 x_2 x_1 x_4$. Relation (2) gets the form 
$x_4 x_3 x_2 x_1 x_3 x_2 = x_2 x_4 x_3 x_2 x_1 x_3$.

Now, we will show that Relation (3) is redundant. We can write Relation (3) also as:
$x_4 x_3 x_2 x_1 x_3 x_2 x_3 = x_3 x_4 x_3 x_2 x_1 x_3 x_2$. By Relation (4), 
the left side becomes $x_3 x_2 x_1 x_4 x_3 x_2 x_3$. Now, we can cancel the leftmost $x_3$, 
to get: $x_2 x_1 x_4 x_3 x_2 x_3 = x_4 x_3 x_2 x_1 x_3 x_2$. By Relation (1), this relation 
is identical to Relation (2), and hence Relation (3) is redundant.

Hence, we get:
\begin{cor}
The singular point presented locally by the equation $(2x+y)(2x-y)(y+x^2)(y-x^2) = 0$ 
has the following induced relations:
$$x_4 x_3 x_2 x_1 = x_1 x_4 x_3 x_2 = x_3 x_2 x_1 x_4 \quad ; \quad x_4 x_3 x_2 x_1 x_3 x_2 = x_2 x_4 x_3 x_2 x_1 x_3$$  
where $\{ x_1,x_2,x_3,x_4 \}$ are the generators of the standard g-base.
\end{cor}

If we delete the generator $x_1$ which corresponds 
to the line $2x+y=0$, we get the following set of relations:
$$x_4 x_3 x_2 = x_3 x_2 x_4 \quad ; \quad x_4 x_3 x_2 x_3 x_2 = x_2 x_4 x_3 x_2 x_3$$  
as expected (see Section \ref{3comps_type2}).
We get similar results if we delete the generator $x_4$ which corresponds 
to the line $2x-y=0$ (see Section \ref{3comps_type1}). 

\medskip

Since the Lefschetz diffeomorphism is a half of the action of the braid monodromy, 
we have that:  
\begin{cor}
The Lefschetz diffeomorphism of the singular point presented locally by 
$(2x+y)(2x-y)(y+x^2)(y-x^2) = 0$ is a counterclockwise full-twist
of the points $y=1$ and $y=-1$, and a counterclockwise half-twist of 
the points $y=2$ and $y=-2$ around the points $y=1,-1$. 
\end{cor}

\subsubsection{Second type}\label{4comps_type5}

\begin{figure}[h]
\epsfysize=2.5cm  
\centerline{\epsfbox{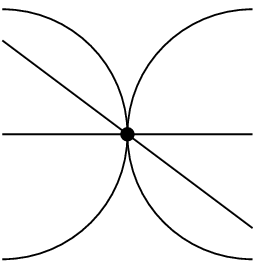}}  
\caption{The singularity of $y(x+2y)(y^2+x)(y^2-x) = 0$ at $(0,0)$}\label{type_i}  
\end{figure}

The local equation of the singularity of the second type is 
$y(x+2y)(y^2+x)(y^2-x) = 0$ (see Figure \ref{type_i}). 
 
There is a major difference 
between this type of singularity and the previous one: In this singularity there are
two ``hidden'' branch points too (as in Section \ref{3comps_type3}). That is, at any fiber, one has four real points 
and two complex points (i.e. complex level $2$), and in each side of the 
singularity, the complex points belong to a different conic (since by the 
singularity, two real points become complex and two complex points become real).   

For computing the braid monodromy in this case, we will use the same trick we
have already used in Section \ref{4comps_type3}:
the two lines in the middle, $x+2y=0$ and $y=0$, can be thought for a moment as a 
``thick'' line which is perpendicular to the tangent direction of the two 
tangented conics.
We have already computed this case (Section \ref{3comps_type3}): we got there that the braid monodromy 
is a $180^{\circ}$ counterclockwise rotation of four points around one fixed point 
at the origin. 
After observing this, we should add into account that the ``thick'' line 
stands for two lines. Hence, the fixed point at the origin is now decomposed 
into two close points, which are doing a counterclockwise full-twist 
(like a usual node).  To summarize, we have: 

\begin{thm}
The local braid monodromy of the singularity presented locally 
by the equation $y(x+2y)(y^2+x)(y^2-x)=0$ is: The four points correspond to $y=1,-1,i$ and 
$y=-i$ do a $180^{\circ}$ counterclockwise rotation, and the two points in the center 
do a counterclockwise full-twist (see Figure \ref{bm_type_i}).

\begin{figure}[h]
\epsfysize=2.8cm  
\centerline{\epsfbox{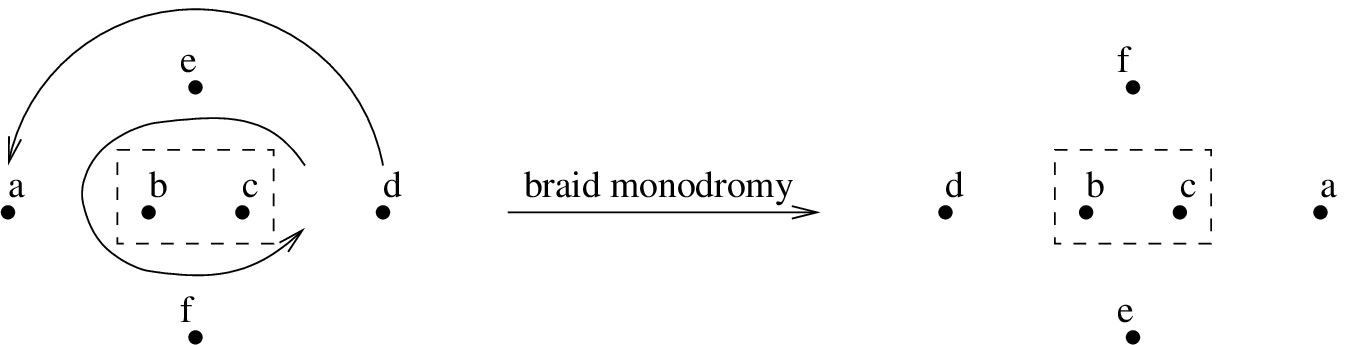}}  
\caption{The action of the local braid monodromy}\label{bm_type_i}  
\end{figure}  

\end{thm}
   
\begin{figure}[h]
\epsfysize=12cm  
\centerline{\epsfbox{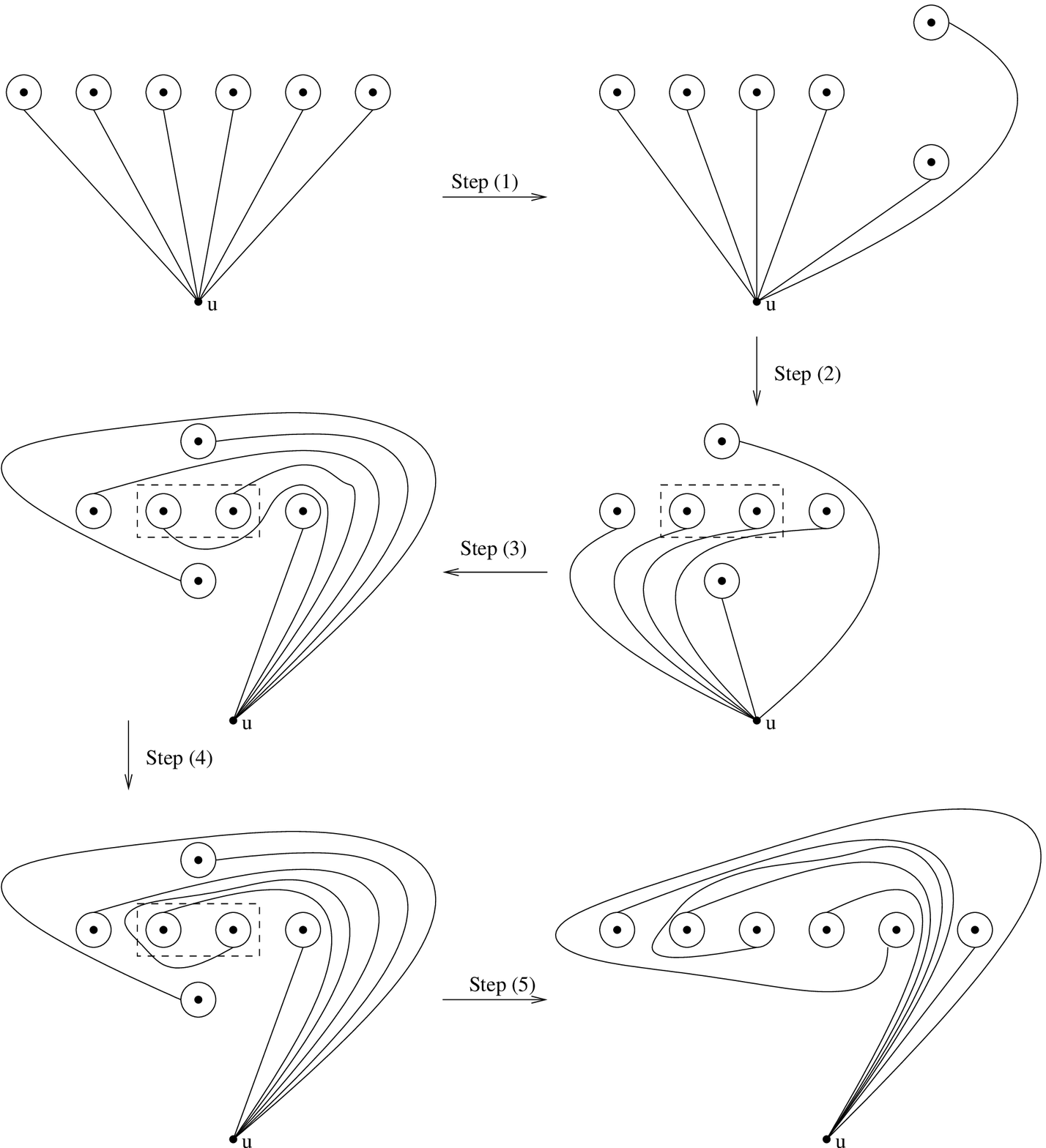}}  
\caption{The g-base obtained from the standard g-base by the action of the local braid monodromy}\label{gbase_type_i}  
\end{figure}

Now, we want to compute the induced relations of this singular point. 
By van Kampen's Theorem, we compute the g-base obtained by applying the action induced
by the local braid monodromy of the singular point on the standard g-base. Since we have 
two complex points in the fiber  before the action of the braid monodromy and after it,
we have to start by rotating the two rightmost points by $90^{\circ}$ 
counterclockwise, for representing the two complex points (see Step (1) in 
Figure \ref{gbase_type_i}). Then, we move the two complex points to be over 
the two middle real points (see Step (2) there). 
Now, we apply the action of the local braid monodromy (Steps (3) and (4)),
and then we return the two new complex points to the right side, and return them to the 
real axis by rotating them by $90^{\circ}$ clockwise (Step (5), see Remark \ref{rem_clockwise}). 

Now, by van Kampen's Theorem, we get the following induced relations from 
the new g-base (where $\{ x_1,x_2,x_3,x_4,x_5,x_6 \}$ are the generators of the standard g-base):
\begin{enumerate}
\item $x_1=x_5 x_4 x_5^{-1}$
\item $x_2=x_5 x_4 x_3 x_2 x_3^{-1} x_4^{-1} x_5^{-1}$
\item $x_3=x_5 x_4 x_3 x_2 x_3 x_2^{-1} x_3^{-1} x_4^{-1} x_5^{-1}$
\item $x_4=x_5 x_4 x_3 x_2 x_1 x_2^{-1} x_3^{-1} x_4^{-1} x_5^{-1}$
\item $x_5=x_6$
\item $x_6=x_6 x_5 x_4 x_3 x_2 x_1 x_5 x_1^{-1} x_2^{-1} x_3^{-1} x_4^{-1} x_5^{-1} x_6^{-1}$ 
\end{enumerate}

From Relation (2), we get $x_5 x_4 x_3 x_2 = x_2 x_5 x_4 x_3$. Using this relation,
Relation (3) becomes: $x_3=x_2 x_5 x_4 x_3 x_4^{-1} x_5^{-1} x_2^{-1}$, and hence
$x_2 x_5 x_4 x_3 = x_3 x_2 x_5 x_4$. Using Relation (2) again, we get 
$x_5 x_4 x_3 x_2 = x_3 x_2 x_5 x_4$.

By Relation (1), Relation (4) can be written: 
$x_4 x_5 x_4 x_3 x_2 x_5 = x_5 x_4 x_3 x_2 x_5 x_4$     

Now, we will show that Relation (6) is redundant. Since $x_5=x_6$, we can simplify 
Relation (6) to the following form: $x_5 x_4 x_3 x_2 x_1 = x_4 x_3 x_2 x_1 x_5$. 
Now, by Relation (1), we get: 
$x_5 x_4 x_3 x_2 x_5 x_4 x_5^{-1} = x_4 x_3 x_2 x_5 x_4 x_5^{-1} x_5$. 
By some simplifications and Relation (3), we get Relation (4), 
and hence Relation (6) is redundant.   

Therefore, we get the following result: 
\begin{cor}
The singular point presented locally by the equation $y(x+2y)(y^2+x)(y^2-x)=0$ has 
the following set of relations:
\begin{enumerate}
\item $x_5 x_4 x_3 x_2 = x_2 x_5 x_4 x_3 = x_3 x_2 x_5 x_4$
\item $x_4 x_5 x_4 x_3 x_2 x_5 = x_5 x_4 x_3 x_2 x_5 x_4$
\item $x_1 = x_5 x_4 x_5^{-1}$
\item $x_5 = x_6$
\end{enumerate}
where $\{ x_1,x_2,x_3,x_4,x_5,x_6 \}$ are the generators of the standard g-base.
\end{cor}

If we delete the generators $x_2$ or $x_3$ which correspond 
to the lines $y=0$ or $x+2y=0$ respectively, we get the set of relations
as in Section \ref{3comps_type3} as expected.

\medskip

Since the Lefschetz diffeomorphism is a half of the action of the braid monodromy, 
we have that:  
\begin{cor}
The Lefschetz diffeomorphism of the singular point presented locally by 
$y(x+2y)(y^2+x)(y^2-x)=0$ is a counterclockwise $90^{\circ}$ rotation of 
the four points around the center (which consists of two points), and then a counterclockwise half-twist 
of the two points in the center. 
\end{cor}

\section*{Acknowledgments} 

The first and second authors wish to thank the Institute of Mathematics, 
Erlangen (Germany), Institut Fourier, Grenoble (France), Hershel Farkas, 
Ehud De-Shalit and Einstein Institute of Mathematics for hosting their stays.

\end{document}